\documentclass[a4paper]{article}
\usepackage{amssymb,amsthm,amsmath}
\usepackage[numbers,sort&compress]{natbib}
\usepackage{esint}
\usepackage{authblk} 
\usepackage{abstract} 

\def\no{\noindent}
\def\endproof{\hphantom{MM}\hfill\llap{$\square$}\goodbreak}

\newenvironment{pf}{{\noindent\bf Proof.\quad}}{\hphantom{MM}\hfill\llap{$\square$}\goodbreak}
\theoremstyle{plain}
\newtheorem{corollary}{Corollary}[section]
\newtheorem{theorem}{Theorem}[section]
\newtheorem{lemma}{Lemma}[section]
\newtheorem{proposition}{Proposition}[section]

\theoremstyle{definition}
\newtheorem{definition}{Definition}[section]

\numberwithin{equation}{section}
\newcommand{\zk}{( }
\newcommand{\yk}{)}
\newcommand {\T}{$T \in (-\infty , \infty]$}

\newcommand{\dz}{( -\Delta_z )^{1/5}} 
\newcommand{\dy}{( -\Delta_y)^{1/3}} 
\newcommand{\spaced}{\mathbb{R}^{1+3d}_T}

\newcommand\cF{\mathcal{F}}
\title{Global $L^p$ Estimates for Kolmogorov--Fokker--Planck Equations in Non-Divergence Form with Four Distinct Scalings}
\author{Liyuan Suo}
\affil{Beijing International Center for Mathematical Research, Peking University, Beijing 100871, P.R.China.}
\affil{liyuansuo@pku.edu.cn}

\begin{document}
\date{}
\maketitle
\begin{abstract}
We study degenerate Kolmogorov--Fokker--Planck operators with four distinct scalings in nondivergence form, where the coefficients $a^{ij}$ are measurable in time and VMO in space. We prove global $L^p$-estimates for $\nabla_x u$, $\dy u$, and $\dz u$. A key contribution is the establishment of a scaling-invariant Poincar\'e inequality for homogeneous solutions, which is fundamental to our regularity analysis.
\end{abstract}



\vskip .3in
\section{Introduction and main results}\hspace*{\parindent}
In this work, we investigate a class of degenerate Kolmogorov--Fokker--Planck (KFP) operators with four distinct scaling structures, expressed in nondivergence form as
\begin{equation}\label{KFP}
Pu = \partial_{t} u - x \cdot \nabla_y u - y \cdot \nabla_z u - \sum_{i,j=1}^d a^{ij}(X)\partial_{x_ix_j}u,
\end{equation}
where $X = (t,x,y,z) \in \mathbb{R}^{1+3d}_T$ with $\mathbb{R}^{1+3d}_T := (-\infty,T) \times \mathbb{R}^{3d}$ for $T \in (-\infty, +\infty]$. The principal coefficients $(a^{ij})_{i,j=1}^d$ are assumed to be bounded measurable functions satisfying the uniform ellipticity condition
\begin{equation}
\lambda|\xi|^2 \leq a^{ij}(X)\xi_i\xi_j \leq \Lambda|\xi|^2 \quad \text{for a.e. } X \in \mathbb{R}^{1+3d}_T, \ \forall \xi \in \mathbb{R}^d,
\end{equation}
with $0 < \lambda \leq \Lambda < \infty$. We denote by $P_0$ the corresponding operator when $a^{ij}$ depend merely on $t$.

The operator \eqref{KFP} constitutes a special case of the general ultraparabolic operator
\begin{equation}\label{L}
L = \partial_t - \sum_{i,j=1}^N b^{ij}x_i\partial_{x_j} - \sum_{i,j=1}^q a^{ij}(t,x)\partial_{x_ix_j}, \quad q \leq N,
\end{equation}
which reduces to a KFP-type operator when the coefficient matrices $(b^{ij})$ and $(a^{ij})$ satisfy specific structural hypotheses (see, e.g., \cite{bramanti1996lpestimates,hormander1967hypoelliptic}). Such operators model essential dynamics in kinetic theory, plasma physics, and related transport-diffusion processes with degenerate structure. The study of these operators has seen significant developments in recent years due to their fundamental role in both theoretical and applied contexts.
\vskip .1in
Due to the structural similarity between KFP operators and second-order parabolic operators, we expect that certain classical results in parabolic theory can be extended to this class of degenerate operators.
%
The De Giorgi-Nash-Moser iteration method, well known in the theory of elliptic and parabolic equations, has made some progress in divergence form KFP equations. Pascucci and Polidoro \cite{pascucci2004moser} successfully achieved local boundedness for weak solutions with measurable coefficients by adapting Moser's iterative scheme. Lunardi \cite{lunardi1997schauder}, Manfredini \cite{manfredini1997dirichlet} and Francesco et al. \citep{farkas2008class} established Schauder estimates for KFP equations. In the case where $a^{ij}(t,x)$ are merely measurable and essentially bounded, Wang and Zhang \citep{wang2009c,wang2010c,zhang2011calpha} obtained $C^{\alpha}$ regularity for weak solutions of the equation and they obtained a particular form of Poincaré inequality satisfied by non-negative weak sub-solution. In 2017, Golse et al \cite{golse2016harnack} proposed an alternate method to establish the H\"{o}lder regularity. From the above, it is evident that the regularity problems of KFP equations share many similarities with those of elliptic and parabolic equations. 

\vskip .1in
The $W^{2,p}$ theory for Kolmogorov-Fokker-Planck equations represents a fundamental research direction, paralleling the classical theory of parabolic equations. For coefficients $a^{ij}(t,x)$ in the VMO class, Bramanti and Cerutti \cite{bramanti1996lpestimates} established interior $L^p$ estimates for second-order derivatives through fundamental solution representations and Calderón-Zygmund theory, effectively generalizing parabolic $W^{2,p}$ estimates to this setting. Subsequent developments by Manfredini and Polidoro \cite{manfredini1998interior} extended these results to divergence-form KFP operators, while Polidoro and Ragusa \cite{polidoro1998sobolev} obtained a priori estimates in Sobolev-Morrey spaces. Bramanti et al. \cite{bramanti2010global, bramanti2013global} established global $L^p$ estimates and the corresponding weak type $(1,1)$ estimates with uniformly continuous coefficients. A central open question concerns the minimal regularity requirements on coefficients for maintaining $W^{2,p}$ regularity. This problem finds its counterpart in parabolic theory, where Krylov \cite{krylov2007parabolic} made significant progress by introducing VMO$_x$-type conditions. His innovative approach, independent of fundamental solutions, combined pointwise sharp function estimates with the Hardy-Littlewood maximal theorem and Fefferman-Stein inequality to establish global $W^{2,p}$ estimates under weakened temporal regularity assumptions. These developments suggest potential pathways for extending $L^p$ regularity theory to KFP equations with less regular coefficients.

\vskip .1in
It should be noted that the aforementioned results primarily address regularity properties of the second-order principal derivatives, while providing no information about the degenerate spatial directions. The strong degeneracy inherent in these operators introduces substantial technical challenges when analyzing the remaining spatial derivatives. In this direction, Bouchut \cite{bouchut2002hypoelliptic} made significant progress in 2002 by establishing fractional derivative estimates for a particular class of KFP equations. Furthermore, in the special case where the coefficients $a^{ij}$ depend merely on the temporal variable, Golse et al. \cite{golse2016harnack} obtained maximal regularity estimates, demonstrating the intricate relationship between temporal regularity and spatial degeneracy in these operators.
%
\vskip .1in
A significant advancement was achieved by Dong and Yastrzhembskiy \cite{dong2022global} in 2022, who extended Krylov's parabolic regularity theory \cite{krylov2007parabolic} to a class of degenerate Kolmogorov-Fokker-Planck equations. Their work considered the non-local evolution equation:
\begin{equation}
\partial_t u - x \cdot \nabla_y u - \sum_{i,j=1}^d a^{ij}(t,x,y)\partial_{x_ix_j}u + \lambda u = f,
\end{equation}
where $\lambda > 0$ represents a damping parameter. The principal coefficients $a^{ij}$ were assumed to belong to a vanishing mean oscillation space VMO$_{x,y}$ with respect to both the $x$ and $y$ variables, generalizing the classical VMO framework to this degenerate setting. This extension required novel techniques to handle the anisotropic scaling and hypoelliptic structure characteristic of KFP operators. Subsequently, they extended their method to divergence-form KFP equations, see \cite{dong2024global} for complete details.
%
%

\vskip .1in
As established in \cite{polidoro1994class}, KFP operators are left-invariant under certain Lie group. The operator \eqref{KFP} investigated in this work represents a significant class of KFP operators. The associated group operation we consider here  takes the form
\[
(t_0,x_0,y_0,z_0) \circ (t,x,y,z) = \left(t+t_0,x+x_0, y + y_0 - t x_0,z+z_0-ty_0+\frac{t^2}{2}x_0\right),
\]
endowed with the anisotropic scaling transformation
\[
(t,x,y,z) \mapsto (r^2t,rx,r^3y,r^5z), \quad r > 0.
\]
This four-parameter scaling generalizes the three-parameter case $(t,x,y)\mapsto (r^2t,rx,r^3y)$ studied by Dong and Yastrzhembskiy \cite{dong2022global}. The primary objective of this paper is to establish global \textit{a priori} estimates for solutions of \eqref{KFP}, specifically proving regularity control for three distinct derivative types: the second-order horizontal derivatives $\nabla_x^2u$, $\dy u$, and $(-\Delta_z)^{1/5}u$. These results substantially extend the existing theory to operators with more complex scaling behavior. 

\vskip .1in
A key aspect of our method is that we establish a kind of Poincaré inequality for the solutions of the homogeneous equation (see Lemma \ref{my2}):
\[  \Vert u \Vert_{L^{2}(Q_{2})} \leq N(d,\delta) \big(\Vert u \Vert_{L^{2}(Q_1)} + \Vert \nabla_zu \Vert_{L^{2}(Q_2)} + \Vert \nabla^2_xu \Vert_{L^{2}(Q_2)} \big). \]
Here, let us revisit the general form of the Poincaré inequality. Suppose \(u(x)\) is a function on \(\mathbb{R}^{d}\), and \(u\in H^{1}(B_2)\), then we have
\[ \Vert u \Vert_{L^{2}(B_2)} \leq N(d) \big( \Vert u \Vert_{L^{2}(B_1)} + \Vert \nabla_x u \Vert_{L^{2}(B_2)}\big).\]

This above inequality implies that if we have the \(L^2\) norm of the derivative of \(u\) in \(B_2\), we can extend the \(L^2\) norm of \(u\) to a bigger domain. We treat the transport term \(\partial_t - x\cdot \nabla_y\) as a whole and utilize the characteristic lines determined by it to connect the points in small regions with those in larger regions, thereby controlling the \(L^2\) norm of \(u\) over larger regions. And the idea of this inequality derives from the Poincaré type inequality in \citep{wang2009c} where Wang and Zhang made use of it to obtain the H\"{o}lder estimates for a class of ultraparabolic equations with measurable coefficients. This fundamental property stems from the intrinsic geometric structure of Kolmogorov–Fokker–Planck operators, characterized by their anisotropic scaling $(t,x,y,z) \mapsto (r^2t,rx,r^3y,r^5z)$, non-commutative Lie group operation, and underlying hypoelliptic nature. The scaling-invariant function spaces and kernel-free methods employed in our work are direct consequences of this geometric foundation. We conjecture that these regularity results can be extended to general Kolmogorov–Fokker–Planck operators with similar underlying geometric structures. 

\vskip.1in 
Recent developments in this direction include the work of Biagi and Bramanti \cite{biagi2024sobolev}, who established global Sobolev estimates for a broader class of KFP operators under the assumption that the coefficients are VMO with respect to the spatial variables. Their approach fundamentally relies on detailed analysis of the fundamental solution for the constant-coefficient model operator. In contrast, our method differs substantially, following the kernel-free techniques developed by Dong \cite{dong2022global}. This alternative approach allows us to obtain new regularity results, particularly in establishing fractional derivative estimates for both $\dy u$ and $\dz u$. These estimates provide finer control over the solution's behavior in the degenerate directions, complementing the existing Sobolev regularity theory.

\vskip .1in
The article is organized as follows: In the remainder of this section, we introduce necessary notations and assumptions, and present our main result (Theorem \ref{T1}). In Section 2, we analyze the case where the coefficients \(a^{ij}\) depend only on \(t\). Using Fourier transform techniques and Parseval's identity, we establish global \(L^2\) estimates. Additionally, we derive localized \(L^2\) estimates, which allow us to prove that \((P_0 + \lambda) C_{0}^{\infty}(\mathbb{R}^{1+3d})\) is dense in \(L^2(\mathbb{R}^{1+3d})\). This result leads to the existence of solutions to the equation, as stated in Theorem \ref{Thm:Solution}. In Section 3, we separately address the Cauchy problem and the homogeneous problem, obtaining pointwise estimates for the sharp functions of \(\partial_z u\) and \(\partial_x^2 u\). These estimates enable us to extend the global results to \(L^p\) spaces for \(p > 1\) via Hardy-Littlewood and Fefferman-Stein type inequalities. Finally, in Section 4, we employ the method of frozen coefficients, locally averaging \(a^{ij}\) with respect to the spatial variables. Combining the results from Section 3 with the VMO conditions satisfied by \(a^{ij}\), we prove our main result, Theorem \ref{T1}.
 
\vskip .3in
\subsection{Notation and the Main Result}\hspace*{\parindent}

For $r>0$ ,  $x_{0}\in \mathbb{R}^{d}$, we set
 \[ B_{r}(x_{0})=\{ x\in \mathbb{R}^{d}:| x-x_{0} | < r \}, B_{r}=B_{r}(0). \]

For $r,R>0$, $X_{0}\in \mathbb{R}^{1+3d}$, we denote
\begin{align*}
   Q_{r,R}(X_{0}) =&\Big\{X \in \mathbb{R}^{1+3d}:-r^2<t-t_{0}<0, |x-x_{0}| < r, |y-y_{0}+ (t-t_{0})x_{0}| < r^3,\\
   & \qquad \qquad\qquad |z-z_0+(t-t_0)y_0-\tfrac{(t-t_0)^2}{2}x_0| < R^{5} \Big\},
\end{align*}
\begin{align*}
   \tilde{Q}_{r,R}(X_{0}) =&\Big\{X \in \mathbb{R}^{1+3d}:|t-t_{0}|<r^2,|x-x_{0}| < r, |y-y_{0}+(t-t_{0})x_{0}| < r^3,\\
   & \qquad \qquad\qquad |z-z_0+(t-t_0)y_0-\tfrac{(t-t_0)^2}{2}x_0| < R^{5}\Big\}.
\end{align*}

For convenience, we abbreviate
\begin{align*}
  &Q_{r}(X_{0})=Q_{r,r}(X_{0}),\quad Q_{r,R}=Q_{r,R}(0),\quad Q_{r}=Q_{r,r}(0),\\
  &\tilde{Q}_{r}(X_{0})=\tilde{Q}_{r,r}(X_{0}),\quad \tilde{Q}_{r,R}=\tilde{Q}_{r,R}(0),\quad 
 \tilde{Q}_{r}=\tilde{Q}_{r,r}(0).
\end{align*}

For any open set $G \subset \mathbb{R}_T^{1+3d}$, we say $u \in S^{p}(G)$ if $u$ satisfies
\[ u,\,\nabla_{x}u, \nabla_x^{2} u, \,\partial_{t} u-x\cdot \nabla_y u-y\cdot \nabla_zu \in L^{p}(G).\]
We define the $S^{p}(G)$ norm of $u$ as 
\begin{align*}
  \|u\|_{S^{p}(G)}:= &\|u\|_{L^{p}(G)} + \|\nabla_{x} u\|_{L^{p}(G)} + \|\nabla_x^{2} u\|_{L^{p}(G)} \\
  & + \|\partial_{t} u-x\cdot \nabla_y u-y\cdot \nabla_zu\|_{L^{p}(G)}.
\end{align*}

For $s \in (0,1/2)$ and $u \in L^{p}(\mathbb{R}^d)$, $(-\Delta_x)^s u$ is understood in the distributional sense:
\[ \left( (-\Delta_x)^s u, \phi \right) = \left( u, (-\Delta_x)^s \phi \right),\quad \phi \in C_{0}^{\infty}(\mathbb{R}^d).\]
When $u$ is a Lipschitz bounded function on $\mathbb{R}^d$, we have the pointwise formula:
\[ (-\Delta_x)^s u(x) = c_{s,d} \int_{\mathbb{R}^d} \frac{u(x)-u(x-\tilde{x})}{|\tilde{x}|^{d+2s}} d\tilde{x},\]
where $c_{s,d}$ depends on $d$ and $s$. For details, see \cite{stinga2019user}.

For any Lebesgue measurable set $\Omega \subset \mathbb{R}^{1+3d}$ with $|\Omega|<\infty$, we denote 
\[(f)_{\Omega}= \fint_{\Omega}f \,dX = |\Omega|^{-1} \int_{\Omega} f \,dX. \]

\vskip .1in
We now state our assumptions on the coefficients.

[$\mathbf{A_1}$] Assume $a^{ij}(X)$, $i,j=1,\cdots,d$ are bounded measurable functions and for some $\delta\in(0,1)$, 
\[ \delta |\xi|^2 \leq a^{ij}(X)\xi_i\xi_j \leq \delta^{-1} |\xi|^2, \qquad \forall X \in \mathbb{R}^{1+3d},\,\xi \in \mathbb{R}^d.\]

The following assumption on $a^{ij}$ can be seen as a $\text{VMO}_{x,y,z}$ condition.
For any $\theta_{0}>0$, there exists $R_{0}=R_{0}(\theta_{0})>0$ such that for all $X_{0}\in\mathbb{R}^{1+3d}$ and $R\in(0,R_{0}]$,
\[ \text{osc}_{x,y,z}(a,Q_{R}(X_{0})) \leq \theta_{0}, \]
where the oscillation is defined by
\begin{equation*}
\begin{aligned}
\text{osc}_{x,y,z}(a,Q_{R}(X_{0})) = \fint_{t_{0}-R^{2}}^{t_{0}} \fint_{D_{R}(X_{0},t)\times D_{R}(X_{0},t)} &|a(t,x_1,y_1,z_1) - a(t,x_2,y_2,z_2)| \\
&dx_1dy_1dz_1dx_2dy_2dz_2dt,
\end{aligned}
\end{equation*}
with the spatial domain
\begin{equation}\label{D}
D_{R}(X_{0},t) = \left\{(x,y,z) : 
\begin{aligned}
&|x-x_0| < R, \\
&|y-y_{0}+(t-t_{0})x_{0}| < R^3, \\
&|z-z_{0}+(t-t_{0})y_{0}-\tfrac{1}{2}(t-t_{0})^{2}x_0| < R^5
\end{aligned}
\right\}.
\end{equation}

We consider the equation with lower order terms:
\begin{equation}\label{eq:main}
Pu + \vec{b}(X)\cdot\nabla_x u + (c(X)+\lambda)u = f,
\end{equation}
where $\lambda > 0$.

[$\mathbf{A_3}$] The coefficients satisfy $|\vec{b}(X)| + |c(X)| \leq L$ for some constant $L>0$ and for all $X\in\mathbb{R}^{1+3d}$.

\begin{definition}
Let $T\in(-\infty,\infty]$. A function $u\in S^{p}(\mathbb{R}^{1+3d}_{T})$ is called a solution of \eqref{eq:main} if the equation holds in the $L^{p}(\mathbb{R}^{1+3d}_{T})$ sense.
\end{definition}

\begin{theorem}\label{T1}
Let $p\in(1,\infty)$, $T\in(-\infty,\infty]$. Under assumptions [$\mathbf{A_1}$] and [$\mathbf{A_3}$], there exists $\theta_{0}=\theta_{0}(d,\delta,L,p)>0$ such that if $\mathbf{A_2}$ holds with this $\theta_{0}$, then:

(i) For $\lambda>\lambda_{0}(d,\delta,L,p)$, we have the a priori estimate:
\begin{equation}\label{mainestimate}
\begin{aligned}
&\lambda\|u\|_{L^{p}} + \lambda^{1/2}\|\nabla_x u\|_{L^{p}} + \|\nabla_x^2 u\|_{L^{p}} \\
&+ \|(-\Delta_z)^{1/5}u\|_{L^{p}} + \|(-\Delta_y)^{1/3}u\|_{L^{p}} \\
&+ \|\nabla_x(-\Delta_y)^{1/6}u\|_{L^{p}} + \|(\partial_t-x\cdot\nabla_y-y\cdot\nabla_z)u\|_{L^{p}} \\
&\leq N\|Pu + \vec{b}\cdot\nabla_x u + (c+\lambda)u\|_{L^{p}},
\end{aligned}
\end{equation}
where $N=N(d,p,\delta,L)$. Moreover, for any $f\in L^{p}(\mathbb{R}^{1+3d}_T)$, \eqref{eq:main} admits a unique solution $u\in S^{p}(\mathbb{R}^{1+3d}_T)$.

(ii) For the Cauchy problem \eqref{eq:Cauchy} with $f\in L^{p}((S,T)\times\mathbb{R}^{3d})$, there exists a unique solution $u\in S^{p}$ satisfying:
\begin{equation}\label{Cauchyu}
\begin{aligned}
&\|u\|_{L^p} + \|\nabla_x u\|_{L^p} + \|\nabla_x(-\Delta_y)^{1/6}u\|_{L^p} \\
&+ \|\nabla_x^2 u\|_{L^p} + \|(-\Delta_y)^{1/3}u\|_{L^p} \\
&+ \|(-\Delta_z)^{1/5}u\|_{L^p} + \|(\partial_t-x\cdot\nabla_y-y\cdot\nabla_z)u\|_{L^p} \\
&\leq N\|f\|_{L^p},
\end{aligned}
\end{equation}
where $N=N(d,\delta,p,T-S)$.
\end{theorem}
\vskip .1in

The core of our approach relies on the Hardy-Littlewood inequality and the Fefferman-Stein inequality. More precisely, we employ the following fundamental estimates.
\begin{lemma}
  Let $c\geq 1$, $ T\in (-\infty,\infty]$. Suppose $f \in L^{p}(\mathbb{R}^{1+3d}_T)$, then we have
  \begin{enumerate}
    \item \textit{Hardy-Littlewood  }
     \[\Vert \mathcal {M}_{c,T} f \Vert_{L^p(\mathbb{R}^{1+3d}_T)} \leq N(d,p) \Vert f \Vert_{L^p(\mathbb{R}^{1+3d}_T)},\]\quad
     
     where \[ \mathcal {M}_{c,T} f (X_{0}) = \sup_{r>0} \fint_{Q_{r,cr}(X_{0})} |f(X)| dX,\qquad \mathcal {M} _{T}= \mathcal {M}_{1,T}, \]
    \item \textit{Fefferman-Stein }
    \[\Vert f \Vert_{L^p(\mathbb{R}^{1+3d}_T)} \leq N(d,p) \Vert  f^{\sharp}_{T} \Vert_{L^p(\mathbb{R}^{1+3d}_T)},\]
    where
\[ f^{\sharp}_{T}  (X_{0}) =\sup_{r>0} \fint_{Q_{r}( X_{0})} |f(X)-(f)_{Q_r ( X_{0}) }| dX.\]
  \end{enumerate}
\end{lemma}

The proof of the above Lemma can be found in the Theorem $7.11$ of \cite{dong2019fully} or \cite{dong2022global}.

\vskip .1in
Next, we introduce the translation and dilation operations for the equation, which will be frequently employed in subsequent analysis. For a fixed point $(t_0,x_0,y_0,z_0) \in \mathbb{R}^{1+3d}$ and scaling parameter $r > 0$, we define the transformed coordinates:

\begin{equation}
\tilde{X} = \left( t_0 + r^2 t,\ x_0 + r x,\ y_0 + r^3 y - r^2 t x_0,\ z_0 + r^5 z - r^2 t y_0 + \frac{r^4 t^2}{2} x_0 \right).
\end{equation}

Given a solution $u$, we define the rescaled function $\tilde{u}(X) := u(\tilde{X})$. A direct computation yields the transformation law for the operator:

\begin{equation}
\left( \partial_t - x \cdot \nabla_y - y \cdot \nabla_z - \sum_{i,j=1}^d a^{ij}(\tilde{X}) \partial_{x_i x_j} \right) \tilde{u}(X) = r^2 (Pu)(\tilde{X}).
\end{equation}

This scaling property reflects the intrinsic anisotropic structure of the Kolmogorov-Fokker-Planck operator, where the temporal and spatial variables scale differently according to the parabolic nature of the equation.


\vskip .3in
\section{ $S^2$ estimate}\hspace*{\parindent}
In this section, we analyze the special case where the coefficients $a^{ij}$ depend solely on the temporal variable $t$. Applying the Fourier transform with respect to all spatial variables $(x,y,z) \in \mathbb{R}^{3d}$, we reduce the partial differential equation to a first-order evolution equation in Fourier space. This reduction enables us to employ the method of characteristics to derive precise $L^2$ estimates for solutions. The main results of this section are summarized as follows.

\begin{theorem}\label{Thm:L2}
   For any $\lambda \geq 0$, $u \in S^2(\mathbb{R}^{1+3d}_T)$, we have the following estimate
  \begin{equation}\label{L2}
  \begin{aligned}
     &\lambda  \Vert u \Vert _{L^2(\mathbb{R}^{1+3d}_T)} + \lambda^{1/2}  \Vert \nabla_x u \Vert _{L^2(\mathbb{R}^{1+3d}_T)}+\Vert \nabla_x(-\Delta_y)^{1/6} u \Vert _{L^2(\mathbb{R}^{1+3d}_T)}\\
     + & \Vert \nabla^{2}_{x} u \Vert _{L^2(\mathbb{R}^{1+3d}_T)}+ \Vert (-\Delta_y)^{1/3} u \Vert _{L^2(\mathbb{R}^{1+3d}_T)}+ \Vert (-\Delta_z)^{1/5} u \Vert _{L^2(\mathbb{R}^{1+3d}_T)} \\
    +& \Vert (\partial_{t}-x\cdot\nabla_{y}-y\cdot\nabla_{z}) u \Vert  _{L^2(\mathbb{R}^{1+3d}_T)} \leq N(d,\delta) \Vert P_0 u +\lambda u \Vert _{L^2(\mathbb{R}^{1+3d}_T)}.
    \end{aligned}
  \end{equation}
\end{theorem}
\begin{theorem}\label{Thm:Solution}
  For a fixed $\lambda > 0$, $T \in (-\infty , \infty]$ and $ f\in L^{2}(\mathbb{R}^{1+3d}_T)$, then the following equation
  \begin{equation}\label{maineq}
  P_0 u +\lambda u =f
  \end{equation}
  has a unique solution $u\in S^2(\mathbb{R}^{1+3d}_T)$.
\end{theorem}

The exponential multiplier technique combined with Theorem \ref{Thm:L2} yields directly the corresponding estimates for zero initial value Cauchy problem on finite intervals, from which we can effortlessly obtain further regularity properties.
\begin{corollary}\label{Cauchy}
   For given numbers $S<T$ and  suppose $f\in L^{2}((S,T)\times\mathbb{R}^{3d})$, the Cauchy initial value problem 
   \begin{equation}\label{eq:Cauchy}
 \begin{cases}
 & P_0 u (X)= f(X), \qquad  X\in (S,T)\times \mathbb{R}^{3d},\\
 & u(S, x,y,z) =0,\qquad (x,y,z) \in \mathbb{R}^{3d}.
 \end{cases}
 \end{equation}
   has a unique solution $u\in S^2((S,T)\times\mathbb{R}^{3d})$. Besides $u$ satisfies
   \begin{equation}\label{Cauchyu}
  \begin{aligned}
     &\Vert u \Vert _{L^2((S,T)\times\mathbb{R}^{3d})} +\Vert \nabla_x u \Vert _{L^2((S,T)\times\mathbb{R}^{3d})}+\Vert \nabla_x(-\Delta_y)^{1/6} u \Vert _{L^2((S,T)\times\mathbb{R}^{3d})}\\
     +& \Vert \nabla^{2}_{x} u \Vert _{L^2((S,T)\times\mathbb{R}^{3d})}+ \Vert (-\Delta_y)^{1/3} u \Vert _{L^2((S,T)\times\mathbb{R}^{3d})}+\Vert (-\Delta_z)^{1/5} u \Vert _{L^2((S,T)\times\mathbb{R}^{3d})} \\
    +& \Vert (\partial_{t}-x\cdot\nabla_{y}-y\cdot\nabla_{z}) u \Vert  _{L^2((S,T)\times\mathbb{R}^{3d})} \leq N(d,\delta,T-S) \Vert f \Vert _{L^2((S,T)\times\mathbb{R}^{3d})}.
    \end{aligned}
  \end{equation}
\end{corollary}
\begin{pf}
Let $\lambda =1$. By Theorem \ref{Thm:Solution}, there exits a $w\in S^2(\mathbb{R}^{1+3d}_T)$ which meets the equation
$$ P_0 w + w = e^{-t} f  \chi_{\{t: S<t<T\}}. $$
In addition, one has
\begin{equation}\label{w}
  \begin{aligned}
     &\Vert w \Vert _{L^2(\mathbb{R}^{1+3d}_T)} +\Vert \nabla_x w \Vert _{L^2(\mathbb{R}^{1+3d}_T)}+\Vert \nabla_x(-\Delta_y)^{1/6} w \Vert _{L^2(\mathbb{R}^{1+3d}_T)} \\
     &+ \Vert \nabla^{2}_{x} w \Vert _{L^2(\mathbb{R}^{1+3d}_T)}+ \Vert (-\Delta_y)^{1/3} w \Vert _{L^2(\mathbb{R}^{1+3d}_T)}+\Vert (-\Delta_z)^{1/5} w \Vert _{L^2(\mathbb{R}^{1+3d}_T)} \\
    &+\Vert (\partial_{t}-x\cdot\nabla_{y}-y\cdot\nabla_{z}) w\Vert  _{L^2(\mathbb{R}^{1+3d}_T)}\\
     \leq & N(d,\delta)\Vert e^{-t} f\chi_{\{t: S<t<T\}}\Vert _{L^2(\mathbb{R}^{1+3d}_T)}\\
\leq & N(d,\delta,T-S) \Vert f\Vert _{L^2((S,T)\times\mathbb{R}^{3d})}.
    \end{aligned}
  \end{equation}
We notice that $e^{-t}f\chi_{\{t: S<t<T\}}\equiv 0$, when $t\leq S$, by the uniqueness of the equation we get that $w=0$, when $t\leq S$. Denote $u(X)=e^t w(X), S\leq t<T$. By direct calculation we have that $u$ is a solution of equation \eqref{eq:Cauchy}. Besides we can get the estimate \eqref{Cauchyu} from \eqref{w}.
\end{pf}
\vskip .1in
Since the coefficients $a^{ij}$ depend only on the temporal variable $t$, we apply the Fourier transform with respect to the spatial variables $(x,y,z)\in\mathbb{R}^{3d}$. Denoting by $U(t,\xi,\eta,\zeta)$ and $F(t,\xi,\eta,\zeta)$ the Fourier transforms of $u(t,x,y,z)$ and $f(t,x,y,z)$ respectively, we obtain the transformed equation:

\begin{equation}\label{eq:U}
\partial_t U + a^{ij}(t)\xi_i\xi_j U + \eta\cdot\nabla_\xi U + \zeta\cdot\nabla_\eta U + \lambda U = F.
\end{equation}

By carefully analyzing the structure of equation \eqref{eq:U}, we employ the method of characteristics to derive an explicit representation for $U$, which subsequently leads to the desired estimates. Following a similar approach to Lemma 4.1 in \cite{dong2022global}, we establish the following result.

\begin{lemma}\label{Lemm:Fourier}
For $\lambda > 0$ and $T \in (-\infty, \infty]$, let $U \in C_b(\mathbb{R}^{1+3d}_T)$ with $\nabla_\xi U, \nabla_\eta U \in C_b(\mathbb{R}^{1+3d}_T)$, and suppose:
\begin{align*}
    \partial_t U &\in L^\infty\big((-\infty, T); C_b(\mathbb{R}^{3d})\big) \cap L^2(\mathbb{R}^{1+3d}_T), \\
    F &\in L^\infty\big((-\infty, T); C_b(\mathbb{R}^{3d})\big) \cap L^2(\mathbb{R}^{1+3d}_T)
\end{align*}
satisfy equation \eqref{eq:U}. Then we derive:
\begin{equation}
\begin{aligned}
&\lambda\Vert U\Vert_{L^{2}(\mathbb{R}^{1+3d}_T)} +\Vert |\xi| ^{2} U\Vert_{L^{2}(\mathbb{R}^{1+3d}_T)} +\Vert |\eta|^{2/3}U\Vert_{L^{2}(\mathbb{R}^{1+3d}_T)}+\Vert |\zeta| ^{2/5}U\Vert_{L^{2}(\mathbb{R}^{1+3d}_T)}\\
+&\| |\zeta|^{1/5}|\xi| U\|_{L^{2}(\mathbb{R}^{1+3d}_T)}
+\Vert |\eta|^{1/3}|\xi|U\Vert_{L^{2}(\mathbb{R}^{1+3d}_T)}\leq N(d,\delta)\Vert F\Vert_{L^{2}(\mathbb{R}^{1+3d}_T)}.
\end{aligned}
\end{equation}
\end{lemma}
\begin{pf}
Using the method of characteristics, we obtain the following explicit representation for $U$:
\begin{equation}\label{rep:U}
  \begin{aligned}
    &U(t, \xi, \eta, \zeta) \\
    =&\int_{-\infty }^{t} e^{-\lambda (t-t')} \exp\Big(-\int_{t'}^{t} a^{ij}(\tau)\xi_i(\tau)\xi_j(\tau) d\tau\Big)\\
    & \qquad\qquad \times F(t', \xi+(t'-t)\eta +\frac{(t'-t)^2}{2}\zeta,\eta +(t'-t) \zeta ,\zeta) dt'.
    \end{aligned}
\end{equation}

The estimates for $\|U\|_{L^{2}(\mathbb{R}^{1+3d}_T)}$, $\||\zeta|^{2/5}U\|_{L^{2}(\mathbb{R}^{1+3d}_T)}$, and $\||\eta|^{2/3}U\|_{L^{2}(\mathbb{R}^{1+3d}_T)}$ follow directly from arguments parallel to those in Lemma 4.1 of \cite{dong2022global}, and thus we omit their proofs here.

\textbullet\textit{Estimate for $\Vert |\xi |^{2}U\Vert_{L^{2}( \mathbb{R}^{1+3d}_T)}.$}  

First by Cauchy-Schwartz inequality, we obtain
 \[ \Vert |\xi |^{2}U\Vert_{L^{2}(\mathbb{R}^{1+3d}_T)} ^2 \leq \int_{\mathbb{R}^{1+3d}_{T}}I_1(X)I_2(X)dX,\]
 where
 \begin{align*}
   & I_1 (X) = \int_{-\infty}^{t} |\xi|^2 e^{-\frac{\delta (t-t')}{1000}|\xi|^2} dt' \leq \int^{\infty}_{0} e^{-\frac{\delta}{1000}t} dt \leq N(\delta),\\
   & I_2(X)  = \int_{-\infty}^{t} |\xi|^2 e^{-\frac{\delta(t-t')}{1000} \Big(|\xi|^2 +(t-t')^2 |\eta| ^2 +(t-t')^4 |\zeta|^2 \Big) } \\
   & \qquad \qquad \qquad \times F^2(t', \xi+(t'-t)\eta+\frac{(t'-t)^2}{2}\zeta,\eta + (t'-t)\zeta ,\zeta) dt'\\
   & \qquad \leq N \int_{-\infty}^{t}\Big(|\xi|^2+ (t-t')^2|\eta|^2  +(t-t')^4|\zeta|^2\Big)\\
   &\qquad \qquad \times e^{-\frac{\delta(t-t')}{2000}\Big(| \xi|^2 + (t-t')^2|\eta| ^2 +(t-t')^4|\zeta|^2\Big)}F^2( t', \xi,\eta ,\zeta) dt', 
 \end{align*}

Taking the advantage of Fubini Theorem, one has
\begin{align*}
   &\Vert |\xi |^{2}U\Vert_{L^{2}(\mathbb{R}^{1+3d}_T)}^2 \\
    \leq & N(\delta) \int_{\mathbb{R}^{1+3d}_T} \int_{-\infty}^{t} (|\xi|^2+ (t-t')^2|\eta|^2+(t-t')^4|\zeta|^2)\\
    & \qquad \qquad \times e^{-\frac{\delta(t-t')}{2000}(| \xi|^2 + (t-t')^2|\eta| ^2 +(t-t')^4|\zeta|^2)}F^2(t', \xi,\eta ,\zeta ) dt' dX\\
  \leq & N(\delta) \int_{\mathbb{R}^{3d}} \Big(\int_{0}^{\infty} \Big(|\xi|^2+ t^2|\eta|^2 +t^4|\zeta|^2\Big)e^{-\frac{\delta}{2000}(t|\xi|^2+ t^3 |\eta| ^2+t^5 |\zeta|^2 )}dt\Big)\\
  & \qquad \qquad \times \left( \int_{-\infty}^{T} F^2(t, \xi,\eta ,\zeta )dt \right)d\xi d\eta d\zeta\\
  \leq & N(\delta) \Vert F\Vert _{L^{2}(\mathbb{R}^{1+3d}_T)} ^2.
\end{align*}

\textbullet\textit{Estimate for $\Vert |\eta|^{1/3}|\xi|U\Vert_{L^{2}(\mathbb{R}^{1+3d}_T)}$ and $\Vert |\eta|^{1/5}|\xi|U\Vert_{L^{2}(\mathbb{R}^{1+3d}_T)}$ } 

Applying the Cauchy-Schwarz inequality, we interpolate between the norms of $|\xi|^{2}U$, $|\eta|^{2/3}U$, and $|\zeta|^{2/5}U$ to obtain the mixed derivative estimates:
\[
\||\eta|^{1/3}|\xi|U\|_{L^{2}(\mathbb{R}^{1+3d}_T)} \quad \text{and} \quad \||\zeta|^{1/5}|\xi|U\|_{L^{2}(\mathbb{R}^{1+3d}_T)}.
\]

Now we have completed the proof of this lemma.
\end{pf}

\vskip .1in
By applying the properties of the Fourier transform and Parseval's identity, we obtain
\begin{align*}
 & \Vert \nabla_x^2u \Vert _{L^2(\mathbb{R}^{1+3d}_T)}=\Vert |\xi| ^{2} U\Vert_{L^{2}(\mathbb{R}^{1+3d}_T)} ,\\
  &\Vert (-\Delta_y)^{1/3} u \Vert _{L^2(\mathbb{R}^{1+3d}_T)}=\Vert |\eta|^{2/3}U\Vert_{L^{2}(\mathbb{R}^{1+3d}_T)},\\
 & \Vert (-\Delta_z)^{1/5} u \Vert _{L^2(\mathbb{R}^{1+3d}_T)} = \Vert |\zeta| ^{2/5}U\Vert_{L^{2}(\mathbb{R}^{1+3d}_T)}.
\end{align*}

Next, we combine the preceding estimates with Lemma \ref{Lemm:Fourier} to prove Theorem \ref{Thm:L2}.

\no {\bf Proof of Theorem \ref{Thm:L2}\quad}
For any $u \in S^{2}(\spaced)$, we construct via mollification (as in Lemma 4.4 of \cite{dong2022global}) a sequence $\{u_n\}_{n=1}^\infty \subset C^\infty(\spaced) \cap S^2(\spaced)$ satisfying
\[
\lim_{n \to \infty} \| u_n - u \|_{S^2(\spaced)} = 0.
\]
The uniform estimates for $u_n$ then follow from Lemma \ref{Lemm:Fourier} and Parseval identity.
\begin{equation}\label{unn}
  \begin{aligned}
     &\lambda  \Vert u_n \Vert _{L^2(\mathbb{R}^{1+3d}_T)} + \lambda^{1/2}  \Vert \nabla_x u_n \Vert _{L^2(\mathbb{R}^{1+3d}_T)}+ \Vert \nabla^{2}_{x} u_n \Vert _{L^2(\mathbb{R}^{1+3d}_T)}\\
      &+ \Vert (-\Delta_y)^{1/3} u_n \Vert _{L^2(\mathbb{R}^{1+3d}_T)} +  \Vert \nabla_x(-\Delta_y)^{1/6} u_n \Vert _{L^2(\mathbb{R}^{1+3d}_T)} \\
     & + \Vert (-\Delta_z)^{1/5} u_n \Vert _{L^2(\mathbb{R}^{1+3d}_T)}+ \Vert (\partial_{t}-x\cdot\nabla_{y}-y\cdot\nabla_{z}) u_n \Vert  _{L^2(\mathbb{R}^{1+3d}_T)}\\
     \leq & N(d,\delta) \Vert P_0 u_n +\lambda u_n \Vert _{L^2(\mathbb{R}^{1+3d}_T)}.
    \end{aligned}
  \end{equation}
  
 Passing to the limit as $n \to \infty$ in the above inequality, we obtain the corresponding $L^2$ estimates for $u$, $\nabla_x u$, and $\nabla_x^2 u$ in \eqref{L2}. 

To establish estimates for the fractional derivative $(-\Delta_z)^{1/5}u$, we employ a duality argument. For any test function $\phi \in C_0^\infty(\spaced)$, the $L^2$-convergence $\|u_n - u\|_{L^2(\spaced)} \to 0$ implies that
  \begin{align*}
&\langle (-\Delta_z)^{1/5} u,\phi\rangle=\langle u,(-\Delta_z)^{1/5} \phi\rangle\\
= &\lim\limits_{n\rightarrow 0}\langle u_n,(-\Delta_z)^{1/5} \phi\rangle
= \lim\limits_{n\rightarrow 0}\langle (-\Delta_z)^{1/5} u_n,(-\Delta_z)^{1/5} \phi\rangle\\
\leq & \Vert\phi \Vert_{L^{2}(\spaced)} \lim\limits_{n\rightarrow 0}\Vert (-\Delta_z)^{1/5} u_n  \Vert_{L^{2}(\spaced)}.
\end{align*}
Combining with \eqref{unn}, we derive
\begin{align*}
\Vert (-\Delta_z)^{1/5} u \Vert _{L^2(\mathbb{R}^{1+3d}_T)}
&\leq  N(d,\delta) \lim\limits_{n\rightarrow 0} \Vert P_0 u_n +\lambda u_n \Vert _{L^2(\mathbb{R}^{1+3d}_T)}\\
&\leq N(d,\delta) \Vert P_0 u +\lambda u\Vert _{L^2(\mathbb{R}^{1+3d}_T)}.
\end{align*}
Following the same approach, we also establish the corresponding estimates for the fractional derivatives:
\[
\|(-\Delta_y)^{1/3}u\|_{L^2(\spaced)} \quad \text{and} \quad \|\nabla_x(-\Delta_y)^{1/6}u\|_{L^2(\spaced)}.
\]
\endproof

\vskip .1in
We can also derive localized $L^2$ estimates through careful selection of appropriate cutoff functions. The proof of this result follows standard techniques similar to those employed in \cite{dong2022global}, and thus we omit the details here for brevity.

\begin{lemma}\label{local}
Let $\lambda \geq 0$, $0 < r_1 < r_2$, and $0 < R_1 < R_2$. Assume $u \in S^2_{\text{loc}}(\mathbb{R}^{1+3d}_0)$ and $f \in L^2_{\text{loc}}(\mathbb{R}^{1+3d}_0)$ satisfy the equation
\[
P_0 u + \lambda u = f.
\]
Then there exists a constant $N = N(d, \delta)$ such that the following estimates holds.
\begin{equation}\label{Ine2}
  \begin{aligned}
    (\romannumeral 1)\quad & ( r_2 - r_1  )^{-1} \Vert \nabla_x u \Vert_{L^{2}(Q_{r_1, R_1})} +  \Vert \nabla_x ^2 u \Vert_{L^{2}( Q_{r_1, R_1})}\\
    \leq & N(d,\delta)\Big((( r_2 - r_1 )^{-2} + r_2 (R_2 - R_1)^{-3} +  R_2 (R_2 - R_1)^{-5})\Vert u \Vert_{L^{2}(Q_{r_2, R_2})}\\
  & \qquad \qquad+ \Vert f \Vert_{L^{2}(Q_{r_2, R_2})} \Big).
  \end{aligned}
\end{equation}

$(\romannumeral 2)$ Denote $C_r = ( -r^2 ,0)\times B_r \times \mathbb{R}^{d} \times \mathbb{R}^{d}.$ Then we get
\begin{equation} \label{Ine3}
  \begin{aligned}
&(r_2 - r_1  ) ^{-1} \Vert \nabla_x u \Vert_{L^{2}(C_{r})} +  \Vert \nabla_x ^2 u \Vert_{L^{2}(C_{r})} \\
\leq& N(d,\delta)\Big(\Vert f \Vert_{L^{2}(C_{r})} +  (r_2 - r_1)^{-2} \Vert u \Vert_{L^{2}(C_{r})}\Big).
  \end{aligned}
\end{equation}
\end{lemma}
\vskip .1in
Building upon the localized $L^2$ estimates established above, we now ready to prove the existence of solutions to equation \eqref{maineq}. 
\begin{lemma}\label{dense2}
For every $\lambda \geq 0$, the set $(P_0 + \lambda)C_0^\infty(\mathbb{R}^{1+3d})$ is dense in $L^2(\mathbb{R}^{1+3d})$.
\end{lemma}

\begin{pf}
We proceed to prove this lemma by contradiction. Suppose that $(P_0 + \lambda)C_0^\infty(\mathbb{R}^{1+3d})$ is not dense in $L^2(\mathbb{R}^{1+3d})$. Then there exists a nonzero function $u \in L^2(\mathbb{R}^{1+3d})$ such that for all test functions $\psi \in C_0^\infty(\mathbb{R}^{1+3d})$, we have the orthogonality relation:
\begin{equation}\label{psi}
    \int_{\mathbb{R}^{1+3d}} (P_0 + \lambda)\psi(X) u(X) \, dX = 0.
\end{equation}

Let $\rho \in C_{0}^{\infty}(\mathbb{R}^{1+3d})$ be a mollifier satisfying $\int_{\mathbb{R}^{1+3d}} \rho(X) \, dX = 1$. For $\epsilon > 0$, we define the rescaled mollifier
\[
\rho^{\epsilon}(t', x', y', z') := \epsilon^{-(2 + 9d)} \rho\left( \frac{t-t'}{\epsilon^2}, \frac{x-x'}{\epsilon^5}, \frac{y-y'}{\epsilon^3}, \frac{z-z'}{\epsilon} \right).
\]

Let
\[ u^{\epsilon} ( X )=\epsilon ^{-2-9d} \int u( t', x', y', z') \rho (\frac{t-t'}{\epsilon^2}, \frac{x-x'}{\epsilon^5}, \frac{y-y'}{\epsilon^3}, \frac{z-z'}{\epsilon} ) dX'.\]


Substituting $\rho^{\epsilon}$ for $\psi$ in \eqref{psi}, we derive the regularized equation for $u^{\epsilon}$:
\begin{equation}
\left( -\partial_{t} + x\cdot\nabla_{y} + y\cdot\nabla_{z} - a^{ij}(t)\partial_{x_i x_j} + \lambda \right) u^{\epsilon}(X) = h^{\epsilon}(X),
\end{equation}
where the source term is given by
\begin{equation}
h^{\epsilon}(X) = \epsilon^{2} \int u(t-\epsilon^2 t', x-\epsilon x', y-\epsilon^3 y', z-\epsilon^5 z') \left( x'\cdot\nabla_{y'} + y'\cdot\nabla_{z'} \right)\rho(t',x',y',z') \, dX'.
\end{equation}

Performing the change of variables $(t,y) \mapsto (-t,-y)$, we define
\begin{equation}
v^{\epsilon}(t,x,y,z) := u^{\epsilon}(-t,x,-y,z),
\end{equation}
which satisfies the equation
\begin{equation}
(P_0 + \lambda) v^{\epsilon}(X) = \tilde{h}^{\epsilon}(X),
\end{equation}
where $\tilde{h}^{\epsilon}(t,x,y,z) = h^{\epsilon}(-t,x,-y,z)$. 

A crucial estimate follows from the structure of $h^{\epsilon}$:
\begin{equation}
\|\tilde{h}^{\epsilon}\|_{L^{2}(\spaced)} \leq N \epsilon^2 \|u\|_{L^{2}(\mathbb{R}^{1+3d})}.
\end{equation}

Applying the local estimate \eqref{L2} with $r>0$ yields:
\begin{equation}\label{Ine5}
\begin{aligned}
\|\nabla_x u^{\epsilon}\|_{L^{2}(Q_r)} 
&\leq \|\nabla_x v^{\epsilon}\|_{L^{2}(Q_r)} \\
&\leq N(d,\delta)\left( r \|\tilde{h}^{\epsilon}\|_{L^{2}(Q_{2r})} + r^{-1} \|v^{\epsilon}\|_{L^{2}(Q_{2r})} \right) \\
&\leq N(d, \delta)(\epsilon^2 r + r^{-1}) \|u\|_{L^{2}(\mathbb{R}^{1+3d})}.
\end{aligned}
\end{equation}

Taking $\epsilon \to 0$ first, we obtain
\begin{equation}
\|\nabla_x u\|_{L^{2}(Q_r)} \leq N(d, \delta) r^{-1} \|u\|_{L^{2}(\mathbb{R}^{1+3d})}.
\end{equation}
Subsequently letting $r \to \infty$ implies $\nabla_x u \equiv 0$ almost everywhere. Consequently, $u \equiv 0$, which contradicts our initial assumption. This completes the proof of the lemma.
\end{pf}

\vskip .1in 
With the density property now established, we turn to the proof of solution existence as asserted in Theorem~\ref{Thm:Solution}.

\no{\bf Proof of Theorem \ref{Thm:Solution}.}\quad 
We divide the proof into two cases according to the time horizon.

\textbf{Case 1: $T=\infty$.} For fixed $\lambda > 0$ and given source term $f \in L^2(\spaced)$, the density lemma guarantees the existence of approximating functions $\{u_n\}_{n=1}^\infty \subset C_0^\infty(\mathbb{R}^{1+3d})$ satisfying
\[
\lim_{n\to\infty} \| (P_0 + \lambda)u_n - f \|_{L^2(\mathbb{R}^{1+3d})} = 0.
\]

Theorem~\ref{Thm:L2} yields the uniform estimate:
\begin{equation}
\begin{aligned}
\lambda \| u_n \|_{L^2} + \| \nabla_x^2 u_n \|_{L^2} &+ \| (\partial_t - x\cdot\nabla_y - y\cdot\nabla_z) u_n \|_{L^2} \\
&\leq N(d,\delta) \| P_0 u_n + \lambda u_n \|_{L^2} \\
&\leq N(d,\delta) \| f \|_{L^2}.
\end{aligned}
\end{equation}

The uniform boundedness of $\| u_n \|_{S^2(\mathbb{R}^{1+3d})}$ implies the existence of a limit function $u \in S^2(\mathbb{R}^{1+3d})$ with weak convergence:
\[
P_0 u_n + \lambda u_n \rightharpoonup P_0 u + \lambda u \quad \text{in } L^2(\mathbb{R}^{1+3d}).
\]
By uniqueness of weak limits, we conclude $P_0 u + \lambda u = f$, establishing the desired solution.

\textbf{Case 2: $T<\infty$.} Applying Case~1 to the truncated problem
\[
P_0 \tilde{u} + \lambda \tilde{u} = f\chi_{t<T},
\]
we obtain a solution $\tilde{u} \in S^2(\spaced)$. Corollary~\ref{Cauchy} shows $\tilde{u}\equiv 0$ for $t\geq T$ since $f\chi_{t<T} = 0$ in this region. The restriction $u := \tilde{u}\chi_{t<T}$ then satisfies
\[
P_0 u + \lambda u = f \quad \text{on } \mathbb{R}^{1+3d}_T,
\]
completing the existence proof.

The combination of both cases establishes Theorem~\ref{Thm:Solution}.
\endproof

\vskip .3in
\section{$S^p$ estimates}\hspace*{\parindent}
In this section, we continue our analysis under the assumption that the coefficients $a^{ij}$ depend only on the time variable $t$. Our main goal is to extend the $L^2$ a priori estimates from Theorem~\ref{Thm:L2} to the more general $L^p$ case for $p > 1$. The proof strategy involves decomposing the solution $u$ into two parts: the part corresponding to the Cauchy problem with zero initial data and the homogeneous part. Through this decomposition, we establish crucial pointwise estimates for the sharp functions of the fractional derivative $(-\Delta_z)^{1/5}u$ and the second-order spatial derivatives $\nabla_x^2 u$, which will ultimately yield the desired $L^p$ estimates via  Hardy–Littlewood and Fefferman–Stein theorems.

\begin{theorem}\label{Thm3.1} 
For any $\lambda \geq 0$, $p\in(1,\infty)$, we have 
$(\romannumeral 1)$\,  Suppose $u \in S^p(\mathbb{R}^{1+3d}_T)$, then
\begin{equation}
\begin{aligned}
& \lambda  \Vert u \Vert _{L^{p} (\mathbb{R}^{1+3d}_{T})} + \lambda^{1/2}  \Vert \nabla_x u \Vert _{L^{p} (\mathbb{R}^{1+3d}_{T})}\\
&+ \Vert \nabla^{2}_{x} u \Vert _{L^{p} (\mathbb{R}^{1+3d}_{T})} + \Vert( -\Delta_z)^{1/5}u \Vert _{L^{p} (\mathbb{R}^{1+3d}_{T})} + \Vert( -\Delta_y)^{1/3}u \Vert _{L^{p} (\mathbb{R}^{1+3d}_{T})}\\
&+ \Vert\nabla_x( -\Delta_y)^{1/6}u \Vert _{L^{p} (\mathbb{R}^{1+3d}_{T})} + \Vert (\partial_{t}-x\cdot\nabla_{y}-y\cdot\nabla_{z}) u \Vert_{L^{p} (\mathbb{R}^{1+3d}_{T})}\\
\leq & N(d,p,\delta) \Vert P_0 u +\lambda u \Vert _{L^{p} (\mathbb{R}^{1+3d}_{T})}.
\end{aligned}
\end{equation}
  $(\romannumeral 2)$ Suppose $ f\in L^{p}(\mathbb{R}^{1+3d}_T)$, then the equation
  \begin{equation}\label{maineq}
  P_0 u +\lambda u =f
  \end{equation}
  has a unique solution $u\in S^p(\mathbb{R}^{1+3d}_T)$.
\end{theorem}
\vskip .1in
Following the same proof method as Corollary \ref{Cauchy}, but now using part (i) of Theorem \ref{Thm3.1} in place of Theorems \ref{Thm:L2} and \ref{Thm:Solution}, we obtain similar $L^p$ estimates for the Cauchy problem on a finite time interval with zero initial condition.

\begin{corollary}\label{Cauchyp}
   For given numbers $S<T$ and $p \in (1,\infty)$, if $f \in L^{p}((S,T)\times\mathbb{R}^{3d})$, then the Cauchy problem 
   \begin{equation}\label{eq:Cauchyp}
 \begin{cases}
 & P_0 u (X)= f(X), \quad X\in (S,T)\times \mathbb{R}^{3d},\\
 & u(S, x,y,z) =0,\quad (x,y,z) \in \mathbb{R}^{3d}.
 \end{cases}
 \end{equation}
   has exactly one solution $u\in S^p((S,T)\times\mathbb{R}^{3d})$. Also,
   \begin{equation}\label{Cauchyup}
  \begin{aligned}
     &\Vert u \Vert _{L^p} +\Vert \nabla_x u \Vert _{L^2}+\Vert \nabla_x(-\Delta_y)^{1/6} u \Vert _{L^2}\\
     +& \Vert \nabla^{2}_{x} u \Vert _{L^2}+ \Vert (-\Delta_y)^{1/3} u \Vert _{L^2}+\Vert (-\Delta_z)^{1/5} u \Vert _{L^2} \\
    +& \Vert (\partial_{t}-x\cdot\nabla_{y}-y\cdot\nabla_{z}) u \Vert _{L^2} \leq N \Vert f \Vert _{L^2}.
    \end{aligned}
  \end{equation}
\end{corollary}

\subsection{Cauchy problem with zero initial data}
\begin{lemma}\label{F}
Take $R \geq 1$. Let $f \in L^{2} (\mathbb{R}^{1+3d})$ with support in $(-1, 0)\times B_1 \times B_1 \times \mathbb{R}^{d}$. If $u \in S^{2}((-1, 0)\times \mathbb{R}^{3d})$ solves 
\begin{equation}\label{Chapter3:Eq:Cauchy}
\begin{cases}
 P_0 u (X)=f(X),\quad X\in(-1,0)\times \mathbb{R}^{3d},\\
  u(-1,x,y,z)  =0,\quad (x,y,z)\in\mathbb{R}^{3d},
\end{cases}
\end{equation}
then:
\begin{equation}\label{eq3.4}
\begin{aligned}
&\Vert |u | + |\nabla_x u| + |\nabla_x^2 u| \Vert_{L^{2}((-1,0)\times B_R\times B_{R^3}
\times B_{R^5})}\\
\leq &N \sum\limits_{k=0}^{\infty} 2^{-k( k-1) /4}R^{-k} \Vert f \Vert _{L^{2}(Q_{1,2^{k+1}R})},
\end{aligned}
\end{equation}
\begin{equation}
  \big(| ( -\Delta_z)^{1/5}u |^2 \big) _{Q_{1,R}}^{1/2}\leq N R^{-2}\sum\limits_{k=0}^{\infty} 2^{-2k} (f^2) _{Q_{1,2^kR}}^{1/2}.
\end{equation}
\end{lemma} 

\begin{pf}
$\bullet$ \textit{Estimates for $u$, $\nabla_x u$, $\nabla_x^2 u$.}
  
  Following the approach of Lemma 5.2 in \cite{dong2022global}, we modify the argument by introducing a decomposition of $f$ along the $z$-direction:
  \[ f = f_0 + \sum \limits_{k=1}^{\infty} f_k := f \chi _{\{z\in B_{(2R ) ^5}\}} + \sum \limits_{k=1}^{\infty} f \chi _{\{z\in B_{(2^{k+1}R) ^5}\backslash B_{(2^{k}R ) ^5}\}}. \]
  Clearly,
\begin{equation}\label{L1}
  \lim_{n\rightarrow\infty} \sum\limits_{k=0}^{\infty} f_k = f \text{ in }L^{2}.
\end{equation}

For each $f_k$ in the Cauchy problem \eqref{Chapter3:Eq:Cauchy}, Theorem \ref{Thm:Solution} gives exactly one solution $u_k \in S^{2}((-1,0) \times \mathbb{R}^{3d})$. From Corollary \ref{Cauchy}, we get these bounds for $u_k$:
\begin{equation}\label{uk}
\begin{aligned}
&\Vert |u_k | + |\nabla_x u_k| + |\nabla_x^2 u_k| \Vert_{L^{2}(( -1,0) \times \mathbb{R}^{3d} )}\\
\leq & N\Vert f_k \Vert _{L^{2}(( -1,0) \times \mathbb{R}^{3d} )}.
\end{aligned}
\end{equation}
 
 Combining this inequality with the convergence in \eqref{L1}, we conclude that $u_k$, $\nabla_x u_k$, and $\nabla_x^2 u_k$ converge in $L^2((-1,0) \times \mathbb{R}^{3d})$:
\[ \lim_{n\rightarrow\infty} \sum\limits_{k=0}^{n} u_k = u, \,\lim_{n\rightarrow\infty} \sum\limits_{k=0}^{n} \nabla_xu_k = \nabla_xu, \,\lim_{n\rightarrow\infty} \sum\limits_{k=0}^{n} \nabla_x^2u_k = \nabla_x^2u.\]

Next, we construct a sequence of cutoff functions. For each integer $j \geq 0$, let 
\[
\phi_j(x,y,z) \in C_0^\infty \left( B_{2^{j+1}R} \times B_{(2^{j+1}R)^3} \times B_{(2^{j+1}R)^5} \right)
\]
be a smooth function satisfying:
\[
\phi_j \equiv 1 \quad \text{on} \quad B_{2^{j+1/2}R} \times B_{(2^{j+1/2}R)^3} \times B_{(2^{j+1/2}R)^5}.
\]
 
%
 Denote
 \[ u_{k,j} = u_k \phi_j, \qquad k\geq 0, j= 0, 1, \cdots , k-1. \]
 
 The function $u_{k,j}$ satisfies the equation:
\[
P_0 u_{k,j} = u_k P_0 \phi_j + \phi_j f_k - 2 a^{ij}(t) \partial_{x_i} \phi_j \partial_{x_j} u_k.
\]

Since $\phi_j f_k \equiv 0$ by construction, applying Theorem \ref{Thm:L2} yields:
 \begin{equation}
 \begin{aligned}
     &  \Vert |u_{k,j} | + |\nabla_x u_{k,j}| + |\nabla_x^2 u_{k,j}| \Vert_{L^{2}((-1,0)\times \mathbb{R}^{3d})}\\
      \leq & N\Vert |u_k P_0\phi_j | +|\nabla_x u_k||\nabla_x \phi_j | \Vert_{L^{2}(( -1,0) \times \mathbb{R}^{3d})}.
 \end{aligned}
 \end{equation}
 
 Then we have:
 \begin{equation}\label{ukj}
 \begin{aligned}
     & \Vert |u_{k,j} | +|\nabla_x u_{k,j}| + |\nabla_x^2 u_{k,j}|\Vert_{L^{2}((-1,0)\times \mathbb{R}^{3d})}\\
    \leq & N2^{-j} R^{-1} \Vert |u_k | + |\nabla_x u_k|\Vert_{L^{2}\big((-1,0)\times  B_{2^{j+1}R} \times B_{(2^{j+1}R) ^3}\times B_{( 2^{j+1}R )^5}\big)}.
 \end{aligned}
 \end{equation}
 
 Combining \eqref{uk} with \eqref{ukj}, we get
 \begin{equation}
 \begin{aligned}
     &  \Vert |u_k| + |\nabla_x u_{k}| + |\nabla_x^2 u_{k}| \Vert_{L^{2}(( -1,0)\times B_{R} \times B_{R^3}\times B_{R ^5})}\\
       \leq  & N^k2^{-k(k-1)/2} R^{-k} \Vert f_k \Vert_{L^{2}\big(( -1,0)\times \mathbb{R}^{3}\big)}\\
      \leq & N2^{-k(k-1)/4} R^{-k} \Vert f \Vert_{L^{2}(Q_{1, 2^{k+1}R})}.
 \end{aligned}
 \end{equation}
 
Combining the estimate \eqref{uk} for $k=0$ with the triangle inequality yields the desired bound \eqref{eq3.4}.
 
 $\bullet$\textit{Estimate of $( -\Delta _z )^{1/5}u$}.
 
 Consider the equation that $u\phi_0$ satisfies
\[ P_0 ( u\phi_0 ) = f \phi_0 + u P_0\phi_j  - 2 a^{ij}(t)\partial_{x_i} \phi_0 \partial_{x_j} u, \]

From Theorem \ref{Thm:L2} and \eqref{eq3.4}, we have the global estimate for $( -\Delta _z )^{1/5}( u\phi_0 )$
\begin{equation}\label{dz1}
\Vert( -\Delta _z )^{1/5}( u\phi_0 )\Vert_{L^{2}((-1,0)\times \mathbb{R}^{3d}\big)} \leq N\sum\limits_{k=0}^{\infty} 2^{-k(k-1) /4}R^{-k} \Vert f \Vert _{L^{2} ( Q_{1,2^{k+1}R})}.
\end{equation}

Next we consider the commutator to get the local estimate of $( -\Delta _z )^{1/5}u$.
\[ \Vert( -\Delta _z )^{1/5}( u\phi_0 ) - \phi_0 ( -\Delta _z )^{1/5} u \Vert_{L^{2}( Q_{1, R})}. \]

Notice that $\phi_0 =1 $ in $B_{2^{1/2}R} \times B_{(2^{1/2}R) ^3}\times B_{(2^{1/2}R)^5}$. Then for any $X\in Q_{1, R}$ and H\"{o}ler inequality we obtain
\begin{align*}
&| ( -\Delta _z )^{1/5}(u\phi_0) - \phi_0 ( -\Delta _z )^{1/5}u |(X)\\
= & c_d\Big|\int _{\mathbb{R} } \frac{u (t, x, y, z-\tilde{z})\phi_0(x, y, z-\tilde{z})-u (t, x, y, z-\tilde{z})\phi_0(x, y, z)}{|\tilde{z}|^{d+2/5}} d\tilde{z}\Big|\\
\leq & N \int_{|z|\geq ( 2^{5/2}+1 )R^5}  \frac{|u (t, x, y, z-\tilde{z})|}{|\tilde{z}|^{d+2/5}} d\tilde{z}\\
\leq & N \sum\limits_{k=0}^{\infty} \int_{2^{5k}R^{5}\leq |\tilde{z}| \leq2^{5(k+1) }R^{5}} \frac{|u (t, x, y, z-\tilde{z})|}{|\tilde{z}|^{d+2/5}} d\tilde{z}\\
\leq & N  \sum\limits_{k=0}^{\infty} 2^{-\frac{5kd}{2}-2k}R^{-\frac{5d}{2}-2} \big(\int_{2^{5k}R^{5}\leq |\tilde{z}| \leq2^{5(k+1) }R^{5}} |u (t, x, y, z-\tilde{z})|^2d\tilde{z}\big)^{1/2}.
\end{align*}

 And in $Q_{1,R}$ we have
\begin{equation}
\begin{aligned}
&\Vert ( -\Delta _z )^{1/5}( u\phi_0 )- \phi_0 ( -\Delta _z )^{1/5} u \Vert_{L^2(Q_{1,R})}\\
  \leq & N  \sum\limits_{k=0}^{\infty} 2^{-\frac{5kd}{2}-2k}R^{-\frac{5d}{2}-2} \Big(\int_{|z|\leq R^5} \int_{2^{5k}R^{5}\leq |\tilde{z}| \leq2^{5( k+1 )}R^{5}}\\
  &\qquad \qquad \qquad   \Vert u (\cdot,z-\tilde{z}) \Vert _{L^2({( -1,0 ) \times B_1 \times B_1})}d\tilde{z} dz \Big) ^{1/2}\\
 \leq &N \sum\limits_{k=0}^{\infty} 2^{-\frac{5kd}{2}-2k}R^{-2}\big(\int_{|z| \leq2^{5( k+2 )}R^{5}}   \Vert u (\cdot,z) \Vert^2_{L^2 ( ( -1,0 ) \times B_1 \times B_1 )} dz \big) ^{1/2}\\
 \leq &N \sum\limits_{k=0}^{\infty} 2^{-\frac{5kd}{2}-2k}R^{-2}\Vert u \Vert_{L^2( Q_{1,2^k R}) }.
\end{aligned}
\end{equation}

Replacing $R$ with $2^k R$ in \eqref{eq3.4} where we obtain estimates for $\| u \|_{L^2( Q_{1,2^k R})}$ and exchanging the order of summation yields:

\begin{align*}
 & N \sum\limits_{k=0}^{\infty} 2^{-\frac{5kd}{2}-2k}R^{-2} \sum\limits_{l=0}^{\infty} 2^{\frac{-l( l-1)}{4}} ( 2^k R) ^{-l} \Vert f \Vert _{L^2(Q_{1,2^{k+l+1} R})}\\
 \leq & NR^{-2}R^{\frac{5d}{2}} \sum\limits_{k=0}^{\infty} \sum\limits_{l=0}^{\infty}2^{-2k} ( |f|^2 ) ^{1/2}_{Q_{1,2^{k+l+1}R}}\\
 \leq  &NR^{-2}R^{\frac{5d}{2}}  \sum\limits_{l=0}^{\infty} 2^{-2l}(|f|^2) ^{1/2}_{Q_{1,2^{l+1}R}}.
\end{align*}

Combining this result with the estimate \eqref{dz1}, we obtain the required bound for $(-\Delta_z)^{1/5}u$. This completes the proof of the desired estimate.
\end{pf}

\vskip .1in
The preceding lemma establishes local estimates for $\nabla_x^2 u$ and $(-\Delta_z)^{1/5}u$ in the context of the Cauchy problem with zero initial data. Since $(-\Delta_z)^{1/5}u$ is a non-local operator, we must carefully handle its decomposition along the $z$-direction. For solutions of the homogeneous equation $P_0u = 0$, we follow the parabolic approach by first proving interior estimates for higher-order derivatives. Starting from the $L^2$ estimates of $(-\Delta_z)^{1/5}u$ in Theorem~\ref{Thm:L2}, we consider the equation satisfied by $(-\Delta_z)^{1/5}u$ to derive estimates for $(-\Delta_z)^{2/5}u$. Then from the equation for $(-\Delta_z)^{2/5}u$ we obtain control of $(-\Delta_z)^{3/5}u$. The critical threshold $2 \times \frac{3}{5} > 1$ enables us to establish estimates for $\nabla_z u$ through interpolation, and similar arguments apply to derive the corresponding estimates for $\nabla_y u$.

\vskip .1in
\subsection{Homogeneous equation}
\begin{lemma}\label{Lyz}
Suppose $u\in S^{2}_{loc} (\mathbb{R}^{1+3d}_0)$ and 
\[ P_0 u = 0, \quad \text{in }Q_ {1} .\]

Then for  $0 < r < R \leq 1$, we have
\begin{equation}
\Vert \nabla_z u \Vert _{L^{2}(Q_r ) } + \Vert \nabla_y u \Vert _{L^{2}(Q_r )}\leq N(d, \delta,r, R) \Vert u \Vert _{L^{2}(Q_R)}.
\end{equation}
\end{lemma}

\begin{pf}
Choose \( r < r_1 < r_2 < R \) and construct cutoff functions \(\rho \in C_0^\infty((-r_1^2, 0) \times B_{r_1})\) with \(\rho \equiv 1\) on \((-r^2, 0) \times B_r\) and \(\psi \in C_0^\infty(B_{r_1^3} \times B_{r_1^5})\) with \(\psi \equiv 1\) on \(B_{r^3} \times B_{r^5}\). The product \(\phi(X) := \rho(t,x)\psi(y,z)\) then yields a cutoff function supported in \(Q_{r_1}\) that equals 1 on \(Q_r\).

Observe that $u\phi$ satisfies the equation
\begin{equation*}\label{Le1}
  P_0 (u\phi )= u P_0 \phi - 2a^{ij} \nabla_{x_i}u\nabla_{x_j}\phi.
\end{equation*}

$\bullet$\textit{ Estimate of  $\nabla_z u.$ } 

From Theorem \ref{Thm:L2}, for $( -\Delta _z )^{1/5}( u\phi)$ we have
 \begin{equation}\label{name1}
\begin{aligned}
  &\Vert ( -\Delta _z )^{1/5}( u\phi)\Vert_{L^{2}(\mathbb{R}_0^{1+3d})}\\
   \leq & N \Vert  u P_0 \phi  \Vert_{L^{2}(\mathbb{R}_0^{1+3d})} + N \Vert  2a^{ij} \nabla_{x_i}u\nabla_{x_j}\phi\Vert_{L^{2}(\mathbb{R}_0^{1+3d})}.
\end{aligned}
\end{equation}

By \eqref{Ine2} of Lemma \ref{local}, we get
\[ \Vert 2a^{ij} \nabla_{x_i}u\nabla_{x_j}\phi\Vert_{L^{2}(\mathbb{R}_0^{1+3d})} \leq N \Vert u \Vert_{L^2(Q_R)}. \]

Substituting the above estimates into \eqref{name1}, we obtain a global estimate for $( -\Delta _z )^{1/5}( u\phi)$:

\begin{equation}\label{L1.1}
\Vert  ( -\Delta _z )^{1/5}(u\phi)\Vert_{L^{2}(\mathbb{R}_0^{1+3d})} \leq  N \Vert u \Vert_{L^2( Q_R)}.
\end{equation}

Next, we consider the function $\omega_1 := ( -\Delta _z )^{1/5} (u\phi )$. Notice that $P_0 \dz = \dz P_0$, and
\begin{equation*}\label{L1.2}
P_0 \omega_1 = ( -\Delta _z )^{1/5} (uP_0\phi )- 2 a^{ij}\nabla_{x_i}  \rho \nabla_{x_j}  ( -\Delta _z )^{1/5}  ( u \psi).
\end{equation*}

Due to Theorem \ref{Thm:L2}, we get the estimate for $(-\Delta _z )^{1/5}\omega_1 =( -\Delta _z )^{2/5} (u\phi )$ 
\begin{equation}\label{L1.3}
  \begin{aligned}
    &\Vert ( -\Delta _z )^{2/5}(u\phi) \Vert_{L^{2}(\mathbb{R}_0^{1+3d})}\\
     \leq &N \Vert( -\Delta _z )^{1/5} (uP_0\phi)\Vert_{L^{2}(\mathbb{R}_0^{1+3d})} 
    + N \Vert2 a^{ij}\nabla_{x_i}  \rho \nabla_{x_j} ( -\Delta _z )^{1/5}  ( u \psi)\Vert_{L^{2}(\mathbb{R}_0^{1+3d})}.
    \end{aligned}
\end{equation}

Denote
\[ I_1 = \Vert( -\Delta _z )^{1/5} (uP_0\phi)\Vert_{L^{2}(\mathbb{R}_0^{1+3d})}, \]
\[ I_2 = \Vert2 a^{ij}\nabla_{x_i}  \rho \nabla_{x_j} ( -\Delta _z )^{1/5} ( u \psi )\Vert_{L^{2}(\mathbb{R}_0^{1+3d})}.\]

For the term $I_1$, $P_0\phi$ can be seen as a cutoff function. Then by \eqref{L1.1}, we get
\begin{equation}\label{L1.4}
  I_1 \leq N \Vert u \Vert_{L^2( Q_R)}.
\end{equation}

Next we consider the term $I_2$. Note that $ ( -\Delta _z )^{1/5} (u\psi )$ satisfies the equation
\[ P_0 ( -\Delta _z )^{1/5} (u\psi )= - \dz \big((x\cdot\nabla_y + y\cdot\nabla_z) \psi u\big) .\]

By Lemma \ref{L2}, we obtain the localized estimate for $\nabla_x\dz (u \psi)$ 
\begin{equation}
  \begin{aligned}
  I_2 \leq& N \Vert \upsilon \dz ( u \psi ) \Vert_{L^{2}(\mathbb{R}_0^{1+3d})} 
+ N  \Vert \upsilon \dz ( (x\partial_y + y\partial_z)\psi  u) \Vert_{L^{2}(\mathbb{R}_0^{1+3d})},
  \end{aligned}
\end{equation}
where $\upsilon( t, x ) \in C_{0}^{\infty} (( -r_2^2, 0)\times B_{r_2})$ and $\upsilon \equiv 1$ in $(-r_1^2, 0) \times B_{r_1}$.

Together with \eqref{L1.1}, we conclude that
\begin{equation}\label{L1.5}
I_2 \leq N \Vert u \Vert_{L^2 (Q_R )}.
\end{equation}

Combine \eqref{L1.4} with \eqref{L1.5}, one has
\begin{equation}\label{w1.1}
  \Vert (-\Delta _z ) ^{2/5} ( u\phi) \Vert_{L^{2}(\mathbb{R}_0^{1+3d})} \leq N   \Vert u \Vert_{L^2 ( Q_R)}.
\end{equation}

So far, we have obtained the estimate for $(-\Delta_z)^{2/5}(u\phi)$, and since $2 \times \frac{2}{5} < 1$, we still cannot  obtain the estimate for $\nabla_z (u\phi)$ by interpolation inequalities. We simply need to repeat the above steps: considering the equation satisfied by $w_2 := (- \Delta_z)^{2/5}(u\phi)$ and then obtaining the estimate for $(- \Delta_z )^{3/5}u$. 

\begin{equation}\label{w2}
P_0 \omega_2  = (- \Delta_z)^{2/5} (uP_0\phi) - 2 a^{ij}\nabla_{x_i}\rho \nabla_{x_j} (- \Delta_z)^{2/5}(u\psi).
\end{equation}

According to Theorem \ref{Thm:L2}, we have 
\begin{align*}
\Vert (-\Delta _z )^{3/5} (u\phi) \Vert_{L^{2}(\mathbb{R}_0^{1+3d})} \leq & N \Vert(- \Delta_z)^{2/5}(uP_0\phi)\Vert_{L^{2}(\mathbb{R}_0^{1+3d})}\\
& + N \Vert 2 a^{ij}\nabla_{x_i}\rho \nabla_{x_j} (- \Delta_z)^{2/5}(u\psi)\Vert_{L^{2}(\mathbb{R}_0^{1+3d}).}
\end{align*}

Denote
\[ I_3 = \Vert(- \Delta_z)^{2/5}( uP_0\phi)\Vert_{L^{2}(\mathbb{R}_0^{1+3d})}, \]
\[ I_4 = \Vert  a^{ij}\nabla_{x_i}\rho \nabla_{x_j} (- \Delta_z)^{2/5}(u\psi)\Vert_{L^{2}(\mathbb{R}_0^{1+3d})}. \]

Then by \eqref{w1.1}, we have 
\begin{equation}\label{w2.2}
  I_3 \leq N \Vert u \Vert_{L^2({Q_R})}.
\end{equation}

For term $I_4$, the function $(- \Delta_z )^{2/5} (u \psi)$ solves the equation
\[ P_0 (- \Delta_z )^{2/5}(u \psi ) = -
(- \Delta_z )^{2/5} (( x\cdot\nabla_y + y\cdot\nabla_z ) \psi u ).\]

By Lemma \ref{local}, we obtain
\begin{align*}
   I_4 \leq &  N \Vert \upsilon(- \Delta_z )^{2/5}( u \psi ) \Vert_{L^{2}(\mathbb{R}_0^{1+3d})}
    + N  \Vert \upsilon (- \Delta_z )^{2/5} (x\cdot\nabla_y + y\cdot\nabla_z (\psi u ) \Vert_{L^{2}(\mathbb{R}_0^{1+3d})}.
\end{align*}

Again by \eqref{w1.1}, we obtian
\begin{equation}\label{w2.3}
I_4 \leq N \Vert u \Vert_{L^2 (Q_R)}.
\end{equation}

Combine $I_3$ with $I_4$, now we conclude that 
\begin{equation}\label{w2.4}
  \Vert (-\Delta _z ) ^{3/5} (u\phi) \Vert_{L^{2}(\mathbb{R}_0^{1+3d})} \leq N   \Vert u \Vert_{L^2 (Q_R)}. 
\end{equation}   

Using \eqref{w2.4} and interpolation inequality, 
one has
\begin{align*}
&\Vert (1- \Delta_z )^{3/5}( u\phi) \Vert_{L^{2}(\mathbb{R}_0^{1+3d})}\\
 \leq & N \Vert u\phi\Vert_{L^{2}(\mathbb{R}_0^{1+3d})} + N \Vert (- \Delta_z)^{3/5}( u\phi) \Vert_{L^{2}(\mathbb{R}_{0}^{1+3d})} \\
\leq & N \Vert u \Vert_{L^2 (Q_R)}.
\end{align*}

Then we obtain 
\begin{align*}
  \Vert \nabla_z u \Vert_{L^2(Q_r) }   \leq &  \Vert \nabla_z (u\phi ) \Vert_{L^2(\mathbb{R}_0^{1+3d})} 
    \leq  N  \Vert (1- \Delta_z ) ^{3/5}( u\phi) \Vert_{L^{2}(\mathbb{R}_0^{1+3d})}\\
    \leq &N \Vert u \Vert_{L^2 (Q_R)} .
\end{align*}

$\bullet$\textit{Estimate of  $\nabla_y u$.} 
Next, we adapt the same method to estimate $\nabla_y u$. Observing the commutation relation
\[
P_0 \nabla_y = \nabla_y P_0 + [\nabla_y, y] \cdot \nabla_z,
\]
we note the appearance of an additional commutator term $[\nabla_y, y] \cdot \nabla_z$. This term must be treated carefully in our estimates. Crucially, since it involves $\nabla_z$, we can control it using the previously established estimates for $\nabla_z u$.The proof follows similar arguments, though we must account for the additional commutator term $[\nabla_y, y]\cdot\nabla_z$ that arises when applying $P_0$ to $\nabla_y u$. This term is controlled using our previous estimates for $\nabla_z u$.

Parallel to the estimate for $\dz (u\phi)$, we also obtain the corresponding bound for $\dy (u\phi)$:
\begin{equation}\label{wy1}
\Vert \dy(u\phi)\Vert_{L^{2}(\mathbb{R}_0^{1+3d})} \leq  N \Vert u \Vert_{L^2(Q_R)}.
\end{equation}

Furthermore, the function $\omega_3 = \dz(u\phi)$ meets the equation
\begin{align*}\label{wy2}
P_0 \omega_3 =& \dy (uP_0\phi)-  2a^{ij} \partial_{x_i}\rho\partial_{x_j}\dy(u \psi)\\
 & + [ \dy y - y\dy ]\cdot \nabla_z( u\phi).
\end{align*}

Due to Theorem \ref{Thm:L2}, 
\begin{equation}\label{wy3}
\begin{aligned}
&\Vert (-\Delta _y)^{2/3}(u\phi) \Vert_{L^{2}(\mathbb{R}_0^{1+3d})}\\
  \leq & N\Vert\dy (uP_0\phi)\Vert_{L^{2}(\mathbb{R}_0^{1+3d})} 
+ N \Vert2a^{ij} \partial_{x_i}\rho\partial_{x_j}\dy(u \psi)\Vert_{L^{2}(\mathbb{R}_0^{1+3d})} \\
 &+ N \Vert [ \dy y - y\dy ] \cdot\nabla_z(u\phi)\Vert_{L^{2}(\mathbb{R}_0^{1+3d})}.
\end{aligned}
\end{equation}

Denote 
\[ I_ 5= \Vert\dy(uP_0\phi)\Vert_{L^{2}(\mathbb{R}_0^{1+3d})}, \]
\[ I_6 = \Vert [ \dy y - y\dy ] \partial_z(u\phi)\Vert_{L^{2}(\mathbb{R}_0^{1+3d})}, \]
\[ I_7 =\Vert a^{ij} \partial_{x_i}\rho\nabla_{x_j}\dy(u\psi)  \Vert_{L^{2}(\mathbb{R}_0^{1+3d})}.\]

By \eqref{wy1}, we have 
\begin{equation}\label{wy4}
  I_5 \leq N \Vert u \Vert_{L^2(Q_R)}.
\end{equation}

Next we consider the term $I_6$. 
\begin{align*}
  & [ \dy y - y\dy ] \cdot\nabla_z(u\phi)  \\
  = &  \int_{\mathbb{R}^d} \frac{\nabla_{z}(u\phi)(y-\tilde{y})}{\tilde{y}^{2/3}}\,d{\tilde{y}}
\end{align*}
By Young's inequality, we obtain
\begin{equation}
\begin{aligned}
 &\Vert [ \dy y - y\dy ] \cdot\nabla_z(u\phi)  \Vert_{L^{2}(\mathbb{R}^{d})} \\
\leq &  \Vert \nabla_z(u\phi) \Vert_{L^{q}(\mathbb{R}^{d})}
\leq   \Vert \nabla_z(u\phi) \Vert_{L^{2}(Q_{r_2})}.
\end{aligned}
\end{equation}
where $\frac{1}{q}=\frac{1}{2}+1-\frac{2s}{d}>\frac{1}{2}$. Then we conclude that 
\begin{equation}\label{wy5}
I_6 \leq N \Vert u \Vert_{L^{2}(Q_R)}.
\end{equation}

%
%
%
%
Similarly, by examining the governing equation for $\dy(u\psi)$, we also obtain:
\begin{equation}\label{w1}
  \Vert (-\Delta _y)^{2/3}(u\phi)\Vert_{L^{2}(\mathbb{R}_0^{1+3d})} \leq N   \Vert u \Vert_{L^2(Q_R)}.
\end{equation}


At last, we conclude that 
\begin{align*}
  \Vert \nabla_y u \Vert_{L^2( Q_r )}  & \leq  \Vert \nabla_y (u\phi)\Vert_{L^2(\mathbb{R}_0^{3d+1})} 
    \leq   N  \Vert (1- \Delta_y) ^{2/3}(u\phi) \Vert_{L^{2}(\mathbb{R}_0^{3d+1})}\\
    &\leq N \Vert u \Vert_{L^2 (Q_R)}.
\end{align*}

Now the Lemma has been proved.
\end{pf}

\vskip .1in
In fact, similar to the homogeneous parabolic equation, we can also obtain interior estimates for higher-order derivatives of $u$ satisfying $P_0 u = 0$, thus deducing the interior continuity of $u$. By induction, we can derive the following lemma.

\begin{lemma}\label{smooth}
For $R\in (1/2, 1)$, $u\in S^2_{loc}(\mathbb{R}_0^{1+3d})$. Suppose $P_0 u = 0 $ in $Q_1$. Then for integers $k, l, m $, we have the following interior estimate
\begin{equation}
  \sup\limits_{Q_{1/2}}|\nabla^m_x\nabla^l_y\nabla^k_z u| +  \sup\limits_{Q_{1/2}}|\partial_t\nabla^m_x\nabla^l_y\nabla^k_z u|  \leq N(d, \delta, R) \Vert u \Vert_{L^2(Q_R)}.
  \end{equation}
  \end{lemma}
\begin{pf}
Noting that $P_0$ commutes with differentiation in the $z$-direction, we derive estimates for arbitrary $z$-derivatives in terms of $\| u \|_{L^2(Q_R)}$. Applying the approach of Lemma~\ref{Lyz}, we then extend these estimates to higher-order derivatives in both the $y$- and $x$-directions. The detailed proof is omitted here for brevity.
\end{pf}
\vskip.1in
Moreover, the control of $\nabla_z u$ requires estimates for $\dz u$, following arguments analogous to Lemma 5.5 in \cite{dong2022global}.

\begin{lemma}\label{z1/5}
Let $r \in (0,1)$, suppose $u\in S^2_{loc}(\mathbb{R}_0^{1+3d})$, and denote $f = P_0 u $. Assume $f =0 $ in $(-1, 0)\times B_1 \times B_1 \times \mathbb{R}^{d}$. Then we have
\begin{equation}\label{dz1/5}
  \Vert \nabla_z u \Vert_{L^2(Q_r)} \leq N( d,\delta, r )\sum\limits_{k=0}^{\infty} 2^{-3k} (|\dz u|^2 )^{1/2}_{Q_{1, 2^k}}.
\end{equation}
\end{lemma}

\begin{pf}
Choose \( r < R < 1 \) and select a cutoff function \(\phi \in C_0^\infty(Q_R)\) and \(\phi \equiv 1\) on \(Q_r\). Throughout the proof, the constant \(N\) may change from line to line but depends solely on \(d\), \(\delta\), \(r\), and \(R\). 

We decompose \(u\) via the Riesz transform \(\mathcal{R}_z\) in the \(z\)-direction, noting the identity \(\mathcal{R}_z (-\Delta_z)^{1/2} = \nabla_z\) that connects fractional and standard derivatives.

We decompose $\nabla_zu$ as follow:
\begin{align*}
 \phi^2 \nabla_zu &= \phi^2 \mathcal{R}_z (-\Delta _z ) ^{1/2} u 
  = \phi^2 \mathcal{R}_z(-\Delta _z )^{3/10} \omega\\
  &= \phi(L\omega + \text{Comm }\omega),
 \end{align*}
where\[\omega = \dz u,\]
\[
 L\omega=\mathcal{R}_z(-\Delta_z)^{3/10}(\phi \omega),
 \]
\[ \text { Comm } \omega =\phi\mathcal{R}_z( -\Delta_z)^{3/10} \omega-\mathcal{R}_z(-\Delta_z)^{3/10}(\phi \omega). \]

$\bullet$\textit{Estimate of $L\omega$.}
In fact, by utilizing the properties of the Riesz transform operator, which maps $L^2$ functions to $L^2$, we have
\begin{align*}
\Vert L\omega  \Vert_{L^2(Q_R) }& \leq \Vert L\omega  \Vert_{L^{2}(\mathbb{R}^{1+3d})}
\leq N \Vert (-\Delta_z) ^{3/10}(\phi \omega)\Vert_{L^{2}(\mathbb{R}^{1+3d}).}
\end{align*}
 
 Notice 
\[ P_0 \omega = 0\quad (-1, 0)  \times B_1 \times B_1 \times \mathbb{R}^{d}. \]

Because $3/10 < 2/5$, the estimation for $(- \Delta_z)^{3/10} (\phi \omega)$ can be obtained similarly to the estimation for $\nabla_z u$ in Lemma \ref{Lyz}. By employing interpolation inequalities, we derive
\begin{equation}
\begin{aligned}
&\Vert (-\Delta_z) ^{3/10}( \phi \omega)\Vert_{L^{2}(\mathbb{R}^{1+3d})}\\
\leq & N \Vert (-\Delta_z) ^{2/5}( \phi \omega) \Vert_{L^{2}(\mathbb{R}^{1+3d})} + \Vert \phi \omega \Vert_{L^{2}(\mathbb{R}^{1+3d})}\\
\leq & N  \Vert \omega \Vert_{L^{2}(Q_R)}. 
\end{aligned}
\end{equation}

Now we get 
\begin{equation}
\Vert L\omega  \Vert_{L^2(Q_R)} \leq N \Vert \omega \Vert_{L^{2}(Q_R)}.
\end{equation}

$\bullet$\textit{Estimate of $\rm{Comm }\, \omega$. }

Next, we utilize the properties of the Riesz transform to estimate $\text{Comm} \, \omega$. Denote
\[ A = \mathcal{R}_z  (-\Delta _z )^{3/10}=  \nabla_z (-\Delta _z) ^{-1/5}.\]

Then we rewrite $\text{Comm}\, \omega$ as
\[ \text{Comm }\omega = \phi A \omega - A (\phi \omega),\]

From the preceding equation, we observe that $\text{Comm} \, \omega$ represents the commutator of $\phi$ with the operator $A$. To proceed, we employ the negative exponential representation of the Riesz potential (as introduced in Definition 1.2 of \cite{grafakos2008classical}), which allows us to express $A$ as a convolution operator. Specifically, for any locally integrable function $\psi \in L^{1}_{\text{loc}}(\mathbb{R}^d)$, we derive the following representation
\begin{equation*}
  (-\Delta _z )^{-1/5} \psi (z)  = c \int_{\mathbb{R}^d} \frac{\psi(\tilde{z})}{|z-\tilde{z}|^{d-2/5}}d\tilde{z}.
\end{equation*}
Then we have 
\begin{equation*}
  \nabla_z(-\Delta _z)^{-1/5} \psi (z)  = c \int_{\mathbb{R}^d} \frac{\psi(\tilde{z})(z-\tilde{z})}{|z-\tilde{z}|^{d-2/5+2}}d\tilde{z}.
\end{equation*}

Thus, for $\text{Comm} \, \omega$, we obtain
\begin{equation}
\begin{aligned}
&|\text{Comm }\omega(X)|\\
\leq & N \int_{\mathbb{R}}\frac{|\omega (t, x, y, z-\tilde{z})||\phi (t, x, y, z) -\phi (t, x, y, z-\tilde{z})| }{|\tilde{z}|^{d+3/5}}\mathrm{d} \tilde{z}\\
= &( \int_{|\tilde{z}| < 2} + \int_{|\tilde{z}| \geq 2}) \frac{|\omega (t, x, y, z-\tilde{z})||\phi (t, x, y, z) -\phi (t, x, y, z-\tilde{z})| }{|\tilde{z}|^{d+3/5} }\mathrm{d} \tilde{z}\\
=&: I_1 (X) + I_2 (X) .
\end{aligned}
\end{equation}

For the term $I_1(X)$, we eliminate the singularity of $|\tilde{z}|^{d+3/5}$ at the origin using the mean value theorem, 
\[ I_1( X )   \leq N \int_{|\tilde{z}| < 2} \frac{|\omega (t, x, y, z-\tilde{z})| } {|\tilde{z}|^{d-2/5}}\mathrm{d} \tilde{z}\]

By Minkowski inequality
\begin{equation}\label{I1}
  \begin{aligned}
    \Vert I_1 \Vert_{L^{2}(Q_R)} & \leq N \int_{|\tilde{z}| <2}\frac{\Vert \omega (\cdot, \cdot-\tilde{z}) \Vert_{L^{2}(Q_R)}}{|\tilde{z}|^{d-2/5}}\mathrm{d} \tilde{z}\\
\leq & N \Vert \omega \Vert_{L^{2}(Q_{1,2})} \int_{|\tilde{z}| < 2}|\tilde{z}|^{-d+2/5}\mathrm{d} \tilde{z} \leq N  \Vert \omega \Vert_{L^{2}(Q_{1,2})}.
\end{aligned}
\end{equation}

Next, let us consider $I_2(X)$. Note that for $X \in Q_R$, $\phi (t,x,y,z-\tilde{z})=0$, then we conclude that

\begin{equation}
\begin{aligned}
I_2( X )  & \leq  N |\phi ( X ) |\int_{|\tilde{z}| \geq 2} \frac{|\omega (t, x, y, z-\tilde{z})|}{|\tilde{z}|^{d+3/5}} \mathrm{d} \tilde{z}
  \leq \sum\limits_{k=0}^{\infty} \int_{2^{5k}\leq |\tilde{z}| < 2^{5(k+1)}}\frac{|\omega (t, x, y, z-\tilde{z})|}{|\tilde{z}|^{d+3/5} }\mathrm{d} \tilde{z}\\
 & \leq \sum\limits_{k=0}^{\infty}  2^{-5kd/2-3k}\Big(\int_{2^{5k}\leq |\tilde{z}| < 2^{5(k+1)}} |\omega (t, x, y, z-\tilde{z})|^2 \mathrm{d} \tilde{z} \Big)^{1/2}.
\end{aligned}
\end{equation}

Then 
\begin{equation}\label{I2}
\begin{aligned}
&\Vert I_2 \Vert_{L^{2}(Q_R)} \\
\leq & N \sum\limits_{k=0}^{\infty} 2^{-5/2kd-3k} \Big(\int_{|z|\leq R^5}\int_{2^{5k}\leq |\tilde{z}| < 2^{5(k+1)}} \Vert\omega (\cdot,z-\tilde{z})\Vert^2_{L^2(-1,0)\times B_1 \times B_1}\mathrm{d} \tilde{z}\Big) ^{1/2}\\
\leq & N\sum\limits_{k=0}^{\infty}2^{-3k}R^{2/5} (|\omega|^2)^{1/2}_{Q_{1, 2^{5k}}}.
\end{aligned}
\end{equation}

Combing \eqref{I1} with \eqref{I2}, we get the desired estimate \eqref{dz1/5}. 
\end{pf}

Next, we establish a crucial inequality analogous to the classical Poincaré inequality for $u$ satisfying the homogeneous equation $P_0u = 0$. We analyze the transport operator \(\partial_t - x\cdot \nabla_y\) as a unified differential structure. By following the characteristic curves generated by this operator, we establish pointwise connections between small and large spatial domains. This approach enables us to bound the \(L^2\) norm of \(u\) over extended regions through local estimates.
\begin{lemma}\label{my2}
Assume $u \in S^{2}(Q_{2})$ and 
\begin{equation}
   \qquad  P_0 u =0 \quad \text{in } \, Q_{2} .
\end{equation}
 
 Then there exists a constant $N=N(d,\delta)$, such that
 \begin{equation}\label{main}
  \Vert u \Vert_{L^{2}(Q_{2})} \leq N(d,\delta) \Big(\Vert u \Vert_{L^{2}(Q_1)} + \Vert \nabla_zu \Vert_{L^{2}(Q_2)} + \Vert \nabla^2_xu \Vert_{L^{2}(Q_2)} \Big).
 \end{equation}
 \end{lemma}
\begin{pf}
First, let us state the general form of the Poincaré inequality. Suppose $u(x)$ is a function on $\mathbb{R}^{d}$, and $u \in H^{1}(B_2)$. Then we have
\begin{equation}\label{p}
  \Vert u \Vert_{L^{2}(B_2)} \leq N(d) \big( \Vert u \Vert_{L^{2}(B_1)} + \Vert \nabla_x u \Vert_{L^{2}(B_2)}\big).
\end{equation}
The proof of this inequality is relatively straightforward, here we omit its proof.

With the help of  Poincaré inequality, we expand the $z$ direction by the boundedness of $\Vert \nabla_zu \Vert_{L^{2}( Q_{2})}$
\begin{equation}
\Vert u \Vert_{L^{2}((-1,0)\times B_1\times B_1\times B_{2^5})} \leq N(d) \Big(\Vert u \Vert_{L^{2}(Q_1)} + \Vert \nabla_z u \Vert_{L^{2}(Q_2)}\Big).
\end{equation}

Next, similarly, we use $\| \nabla^2_xu \|_{L^{2}(Q_2)}$ to expand in the $x$ direction. Firstly we need to obtain an estimate for $\nabla_x u$. Using interpolation inequalities, we obtain
\[ \Vert \nabla_x u \Vert_{L^{2}((-1,0)\times B_1\times B_1\times B_{2^5})} \leq N(d)\big(\Vert u \Vert_{L^{2}((-1,0)\times B_1\times B_1\times B_{2^5})} + \Vert \nabla_x^2 u \Vert_{L^{2}(Q_2)}\big).\]

Then we conclude that
\begin{equation}
\begin{aligned}
 \Vert \nabla_x u \Vert_{L^{2}((-1,0)\times B_2\times B_1\times B_{2^5})}& \leq N(d)\big(\Vert \nabla_xu \Vert_{L^{2}((-1,0)\times B_1\times B_1\times B_{2^5})} + \Vert \nabla_x^2 u \Vert_{L^{2}(Q_2)}\big)\\
 & \leq  N(d) \big(\Vert u \Vert_{L^{2}(Q_1)} + \Vert \nabla_zu \Vert_{L^{2}(Q_2)} + \Vert \nabla^2_xu \Vert_{L^{2}(Q_2)} \big).\
 \end{aligned}
 \end{equation}
 
 Now we expand $x$ direction
 \begin{equation}\label{x}
\begin{aligned}
 &\Vert u \Vert_{L^{2}((-1,0)\times B_2\times B_1\times B_{2^5})}\\
 \leq & N(d)\big(\Vert u \Vert_{L^{2}((-1,0)\times B_1\times B_1\times B_{2^5})} + \Vert \nabla_x u \Vert_{L^{2}((-1,0)\times B_2\times B_1\times B_{2^5})}\big)\\
 \leq & N(d) \big(\Vert u \Vert_{L^{2}(Q_1)} + \Vert \nabla_zu \Vert_{L^{2}(Q_2)} + \Vert \nabla^2_xu \Vert_{L^{2}(Q_2)} \big).
 \end{aligned}
 \end{equation}
 
Observe that we have already extended the domain in $(x,z)$ from $B_1 \times B_1$ to $B_2 \times B_{2^5}$ using the control provided by $\| u \|_{L^{2}(Q_1)} + \| \nabla_z u \|_{L^{2}(Q_2)} + \| \nabla_x^2 u \|_{L^{2}(Q_2)}$. We now use the equation's structure to further expand the domain in the $t-$ and $y$-directions.

In fact, $u$ satisfies the equation
\[ \partial_t u - x\cdot\nabla_y u = y\cdot\nabla_z u + a^{ij}(t)\partial_{x_ix_j} u \quad \text{in }\, Q_2.\]
So denote $g:=y\cdot\nabla_z u + a^{ij}(t)\partial_{x_ix_j} u $, then we have
\[\Vert g \Vert_{L^2(Q_2)}\leq N(\delta) \big(\Vert \nabla_zu \Vert_{L^{2}(Q_2)} + \Vert \nabla^2_xu \Vert_{L^{2}(Q_2)} \big).\]

We shall utilize the characteristic lines of $\partial_t - x \cdot \nabla_y$ and employ an iterative method to gradually expand the region in $t$ and $y$. Suppose for $1\leq r <2-\frac{1}{96}$, we have 
\begin{equation}\label{r}
  \begin{aligned}
 &\Vert u \Vert _{L^{2}((-r^2,0)\times B_2\times B_{r^3}\times B_{2^5})} \\
 \leq & N(d,\delta) \Big(\Vert u \Vert_{L^{2}(Q_1)} + \Vert \nabla_zu \Vert_{L^{2}(Q_2)} + \Vert \nabla^2_xu \Vert_{L^{2}(Q_2)} \Big).
  \end{aligned}
 \end{equation}
 Then for $R=r+\frac{1}{96}$, we obtain
 \begin{equation}\label{R}
  \begin{aligned}
 &\Vert u \Vert _{L^{2}((-R^2,0)\times B_2\times B_{R^3}\times B_{2^5})} \\
 \leq & N(d,\delta) \Big(\Vert u \Vert_{L^{2}(Q_1)} + \Vert \nabla_zu \Vert_{L^{2}(Q_2)} + \Vert \nabla^2_xu \Vert_{L^{2}(Q_2)} \Big).
  \end{aligned}
 \end{equation}
 
 For simplicity, we omit the $z$ variable in the following proof.  For $(t,y) \in (-R^2,0) \times B_{R^3}$. Define $\hat{t}$ as a function of $t$ such that $\hat{t} = t$ if $-R^2 < t \leq -\frac{1}{2}$; $\hat{t} = t - \frac{3}{4}$, if $-\frac{1}{2} < t < 0$ . Note that the choice of $\hat{t}$ ensures $(\hat{t} + \frac{1}{4}, \hat{t} + \frac{1}{2}) \subset (-r^2,0)$. Additionally, let $\hat{x}$ be a function of $t$ and $y$ such that $\hat{x} = \frac{y}{R^3}$, if $-R^2 < t \leq -\frac{1}{2}$; $\hat{x} = -\frac{y}{R^3}$, if $-\frac{1}{2} < t < 0$. Since $|y| \leq R^3$, it follows that $B_{1/2}(\hat{x}) \subset B_2$.
  
   For any $\tilde{t} \in (\hat{t}+1/4,\hat{t}+1/2) $, $\tilde{x} \in B_{1/2}(\hat{x})$, due to the choice of $\hat{x}$, we observe that $(\tilde{t}-t)\hat{x}$ is always opposite in direction to $y$. Direct computation yields
   \[ |y-(\tilde{t}-t)\hat{x}| \leq R^3 - |\tilde{t}-t| \leq r^3 - \frac{1}{2}|\tilde{t}-t|. \]
 So we have
 \begin{equation}\label{y-tx}
 y-(\tilde{t}-t)\tilde{x} \in B_{r^3}.
 \end{equation}
 
Next, we connect $(t,\tilde{x},y)$ and $(\tilde{t},\tilde{x},y-(\tilde{t}-t)\tilde{x})$ by characteristic lines.
\begin{equation}
\begin{aligned}
& u(t,\tilde{x},y)- u(\tilde{t},\tilde{x},y-(\tilde{t}-t)\tilde{x})\\
= & -u(s\tilde{t}+(1-s)t,\tilde{x},y-s(\tilde{t}-t)\tilde{x})|_{s=0}^{1}\\
= & - \int_{0}^{1} (\tilde{t} -t)(\partial_t u - \tilde{x}\cdot\nabla_yu)(s\tilde{t}+(1-s)t,\tilde{x},y-s(\tilde{t}-t)\tilde{x}))ds\\
= & -\int_{0}^{1} (\tilde{t} -t)g(s\tilde{t}+(1-s)t,\tilde{x},y-s(\tilde{t}-t)\tilde{x}))ds.
\end{aligned}
\end{equation} 

Taking the $L^2$ integral of the above expression over $\tilde{t} \in (\hat{t}+\frac{1}{4},\hat{t}+\frac{1}{2})$, $(t,y,\tilde{x}) \in (-R^2,0) \times B_{R^3} \times B_{1/2}(\hat{x})$, and utilizing the Minkowski inequality, we can deduce
\begin{equation}
\begin{aligned}
&\int_{-R^2}^{0}dt\int_{B_{R^3}}dy \int_{B_{1/2}(\hat{x}) } |u(t,\tilde{x},y)|^2d\tilde{x} \\
 \leq &N\int_{-R^2}^{0}dt\int_{\hat{t}+1/4}^{\hat{t}+1/2}d\tilde{t}\int_{B_{R^3}}dy \int_{B_{1/2}(\hat{x}) }  |u(\tilde{t},\tilde{x},y-(\tilde{t}-t)\tilde{x})|^2d\tilde{x} \\
+&N \int_{-R^2}^{0}dt\int_{\hat{t}+1/4}^{\hat{t}+1/2}d\tilde{t}\int_{B_{R^3}}dy \int_{B_{1/2}(\hat{x}) } \big(\int_{0}^{1} (\tilde{t} -t)g(s\tilde{t}+(1-s)t,\tilde{x},y-s(\tilde{t}-t)\tilde{x})ds\big)^2 d\tilde{x}.
\end{aligned}
\end{equation}

Through an appropriate change of variables and by interchanging the order of integration, we obtain
\begin{equation}
\begin{aligned}
&\int_{-R^2}^{0}dt\int_{B_{R^3}}dy \int_{B_{1/2}(\hat{x})}|u(t,\tilde{x},y)|^2d\tilde{x}  \\
\leq & N(d,\delta) \Big(\Vert u \Vert_{L^{2}(Q_1)}+\Vert g \Vert_{L^{2}(Q_2)}\Big)\\
\leq & N(d,\delta) \Big(\Vert u \Vert_{L^{2}(Q_1)} + \Vert \nabla_zu \Vert_{L^{2}(Q_2)} + \Vert \nabla^2_xu \Vert_{L^{2}(Q_2)} \Big).
\end{aligned}
\end{equation}

The left-hand side of the above integral $\| u \|_{L^2}$ is only local with respect to $x$, we can utilize the boundedness of $\nabla_x^2 u$ on $Q_2$ and once again apply the Poincaré inequality to obtain \eqref{R}.

In the above process, we successfully expand $(t,y) \in (-r^2,0) \times B_{r^3}$ to $(t,y) \in (-R^2,0) \times B_{R^3}$. Utilizing \eqref{x}, we start from $r=1$ and iteratively proceed to $R=2$, thus we obtain \eqref{main}. At this point, we have completed the proof of this lemma.
\end{pf}

\vskip .1in
Applying the preceding lemma, we establish interior estimates for $\nabla_x^2 u$ in terms of itself over expanded domains and $\dz u$.

\begin{lemma}\label{x2}
Suppose $u\in S^2_{loc} (\mathbb{R}_0^{1+3d})$, $P_0 u = 0$  in $(-1,0)\times B_1\times B_1\times \mathbb{R}^{d}$. Then for any integral $k, l, m $, we get
\begin{equation}\label{dx2}
  \begin{aligned}
    & \sup\limits_{Q_{1/2}}|\nabla^{m+2}_x\nabla^l_y\nabla^k_z u| +  \sup\limits_{Q_{1/2}}|\partial_t\nabla^{m+2}_x\nabla^l_y\nabla^k_z u| \\
       \leq &  N(d, \delta) \Vert \nabla^2_x u \Vert_{L^2( Q_R )} + N(d, \delta) \sum\limits_{k=0}^{\infty} 2^{-3k} ( |\dz u|^2 )^{1/2}_{Q_{1, 2^k}}.
  \end{aligned}
\end{equation}
\end{lemma}

\begin{pf}
Denote
\[ u_1 (X)= u(X)- (u)_{Q_{r}} - A^jx_j - B^j (tx_j + y_j)- C^{jl}(x_iy_j-x_jy_i),\]
where$A^j, B^j,C^{jl}$ $(i=1,\cdots,d$, $1\leq j<l\leq d)$ are determined by 
 \[ \int_{Q_{r}} x_j u_1 = \int_{Q_{r}} y_j u_1 =\int_{Q_{r}} x_jy_l u_1 = 0.\] 

Notice that
 \[ P_0 u_1 = 0.\]

Then by Lemma \ref{smooth}, we conclude that
\begin{equation}\label{no}
 \sup\limits_{Q_{1/2}}|\nabla^{k+2}_x\nabla^l_y\nabla^m_z u| +  \sup\limits_{Q_{1/2}}|\partial_t\nabla^{k+2}_x\nabla^l_y\nabla^m_z u|  \leq N \Vert u_1 \Vert_{L^{2}(Q_{r})}.
 \end{equation}
 
 Now we claim that 
 \begin{equation}\label{Sobolev}
 \Vert u_1 \Vert_{L^{2}(Q_{r})} \leq N \Vert \nabla^2_x u \Vert_{L^{2}( Q_{R})} + N \Vert \nabla_z u \Vert_{L^{2}(Q_{R})}.
 \end{equation}

 We proof the claim by contradiction. Suppose the assertion is false, then there exists a sequence $\{u^n\} \in S^2_{\text{loc}}(\mathbb{R}_0^{1+3d})$ such that $P_0 u^n = 0$ on $Q_1$. Substituting $u$ with $u^n$ in the definition of $u_1$, we obtain the corresponding $u^n_1$, and
 \begin{equation}\label{un}
 \Vert u^n_1 \Vert_{L^{2}(Q_{r})} > n \big( \Vert \partial^2_x u^n \Vert_{L^{2}( Q_{R})} +  \Vert \partial_z u^n \Vert_{L^{2}(Q_{R})}\big).
\end{equation}

We normalize and suppose $\Vert u^n_1 \Vert_{L^{2}(Q_{r})} =1$ . Then by Lemma \ref{my2}, we get
\[ \Vert u^n_1 \Vert_{L^{2}(Q_{R})} \leq N.\]

Furthermore, by Lemma \ref{smooth}, the uniform boundedness of the $L^2$ norm of $\{u^n_1\}$ over $Q_r$, there exist a $v\in S^{2}(Q_r)$, satisfies $P_0 v = 0 $, and $\nabla^2_x v = \nabla_z v =0 $,
\[ u^n_1 \rightarrow v, \quad \text{in}\, L^ {2}( Q_{r}).\]
 
 Besides, we also have 
 \[ \int_{Q_{r}} v = \int_{Q_{r}} x_j v= \int_{Q_{r}} y_j v =\int_{Q_{r}} x_jy_lv= 0.\]
 
 However, Lemma~\ref{my1} implies $v \equiv 0$, and consequently $\| v \|_{L^2(Q_r)} = 0$, which leads to a contradiction. This establishes the desired claim. Combining the Sobolev-type inequality \eqref{Sobolev} with the fractional estimate \eqref{dz1/5} yields \eqref{no}, thereby completing the proof of the lemma.
 \end{pf}
 
 \vskip .1in
 Then with the help of Lemma \ref{smooth} and Lemma \ref{x2}, we shall obtain pointwise estimates for the sharp functions of $\partial_z u$ and $\nabla_z^2 u$.
 
 \begin{proposition}\label{38}
Let $r>0$, $\upsilon \geq 2$, \T, for fixed $ X_0=(t_0,x_0,y_0,z_0) \in \overline{\mathbb{R}_{T}^{1+3d}}$. Suppose $u \in S^{2}_{loc}(\mathbb{R}_T^{1+3d})$, and $P_0 u =0$ in $(t_0-\upsilon^2r^2,t_0)\times B_{\upsilon r }(x_0)\times B_{\upsilon^3 r^3 }(y_0) $, then there exists a constant $N = N(d,\delta)$, such that
\begin{equation}
  \begin{aligned}
  (\romannumeral1)\quad I_1:&=\Big(|\dz u - (\dz u ) _{Q_r( X_0)}|^2\Big)^{1/2}_{Q_r(X_0)} \\
  &\leq N \upsilon^{-1} \Big(|\dz u|^2\Big)^{1/2}_{Q_{\upsilon r}( X_0)},\\
  (\romannumeral2) \quad I_2:&=\Big( |\nabla_x^2 u - ( \nabla_x^2 u ) _{Q_r( X_0)}|^2\Big)^{1/2}_{Q_r( X_0)} \\
  & \leq N\upsilon^{-1} \Big(|\nabla_x^2 u|^2\Big)^{1/2}_{Q_{\upsilon r}(X_0)} +  N\upsilon^{-1} \sum\limits_{k=0}^{\infty} 2^{-3k} \Big(|\dz u|^2\Big)^{1/2}_{Q_{\upsilon r, 2^k\upsilon r}(X_0)}.
  \end{aligned}
  \end{equation}
\end{proposition} 

\begin{pf}
The translation invariance and scaling symmetry of the equation allow us to reduce the general case to the specific situation where $r=1/{\upsilon}$ and $X_0 = 0$.

Since $\dz P_0 = P_0 \dz$, then we get
\begin{equation}
P_0 (\dz u)= 0, \quad \text{in}\, (-1,0 )\times B_1 \times B_1\times\mathbb{R}^{d}.
\end{equation}

Then with the help of Lemma \ref{smooth}, we obtain
\begin{equation}
\begin{aligned}
I_1 &\leq \sup\limits_{X_1,X_2\leq Q_{1/\upsilon}} |\dz u \zk X_1\yk - \dz u \zk X_2 \yk|\\
& \leq N \upsilon^{-1} \sup\limits_{Q_{1/2}}\zk|\nabla_x \dz u| +  |\nabla_y \dz u| + |\nabla_y \dz u| +|\partial_t \dz u|\yk \\
& \leq N \zk \upsilon^{-1}|\dz u|^2\yk^{1/2}_{Q_{1}}.
\end{aligned}
\end{equation}

For term $I_2$, an analogous application of Lemma~\ref{x2} yields
\begin{equation}
\begin{aligned}
I_2 &\leq \sup\limits_{X_1,X_2\leq Q_{1/\upsilon}} |\nabla^2_x u \zk X_1\yk - \nabla^2_x u \zk X_2 \yk|\\
& \leq N \upsilon^{-1} \sup\limits_{Q_{1/2}}\zk|\nabla_x \nabla^2_x u| +  |\nabla_y \nabla^2_x  u| + |\nabla_y\nabla^2_x  u| +|\partial_t \nabla^2_x  u|\yk \\
& \leq N\upsilon^{-1} \zk |\nabla_x^2 u|^2\yk^{1/2}_{Q_{1}} + N \upsilon^{-1} \sum\limits_{k=0}^{\infty} 2^{-3k} \zk |\dz u|^2\yk^{1/2}_{Q_{1, 2^k}}.
\end{aligned}
\end{equation}

Putting them all together, we have completed the proof of this proposition.
\end{pf}
 
\vskip .2in
\subsection{The proof of Theorem \ref{Thm3.1}}

Having established separate estimates for solutions to both the zero initial-value Cauchy problem and the homogeneous equation, we now combine Lemma~\ref{F} with Proposition~\ref{38} to derive the following a priori estimates for solutions $u$ of $P_0 u = f$.

\begin{proposition}\label{sharp}
Let $r > 0$, $\upsilon \geq 2$, $T \in (-\infty, \infty]$, and $X_0 \in \overline{\mathbb{R}^{1+3d}_T}$ be given. For any solution $u \in S^2(\mathbb{R}^{1+3d}_T)$ of the equation $P_0 u = f$ in $\mathbb{R}^{1+3d}_T$, there exists a constant $N = N(d, \delta)$ such that the following estimates hold.
\begin{equation}\label{sharp:dzu}
  \begin{aligned}
    (\romannumeral 1)\quad&\big(|\dz u - (\dz u) _{Q_r(X_0)}|^2\big)^{1/2}_{Q_r(X_0)}\\
    \leq & N \upsilon^{-1} \big( |\dz u|^2\big)^{1/2}_{Q_{\upsilon r}(X_0)} + N \upsilon^{\frac{2+9d}{2}} \sum\limits_{k=0}^{\infty} 2^{-2k} \big(|f|^2\big)_{Q_{\upsilon r,2^k\upsilon r}(X_0)}^{1/2},\\
    (\romannumeral2)\quad&\big(|\nabla_x^2 u - \big(\nabla_x^2 u\big)_{Q_r(X_0)}|^2\big)^{1/2}_{Q_r(X_0)} \\
   \leq & N \upsilon^{-1} ( |\nabla_x^2 u|^2)^{1/2}_{Q_{\upsilon r}(X_0)} + N \upsilon^{-1} \sum\limits_{k=0}^{\infty} 2^{-3k} \big( |\dz u|^2\big)^{1/2}_{Q_{\upsilon r, 2^k\upsilon r}(X_0)}\\
   & + N\upsilon^{\frac{2+9d}{2}}  \sum\limits_{k=0}^{\infty} 2^{-k} \big(|f|^2\big)  _{Q_{\upsilon r,2^k\upsilon r}(X_0)}^{1/2}.
  \end{aligned}
\end{equation}
\end{proposition} 
\begin{pf}
Without loss of generality, it suffices to consider the normalized case where $r=1$ and $X_0=0$.
Let $\psi$ be a cutoff function for $(t,x,y)$ with $\text{supp} \psi \subset (-(2\upsilon)^2, 0) \times B_{2\upsilon} \times B_{(2\upsilon)^3}$ and $\psi \equiv 1$ on $(-\upsilon^2, 0) \times B_\upsilon \times B_{\upsilon^3}$. Then Theorem~\ref{Cauchy} guarantees the existence of a unique solution $g \in S^2\left((-(2\upsilon)^2, 0) \times \mathbb{R}^{3d}\right)$ to the Cauchy problem:
\begin{equation}
\begin{cases}
P_0 g &= f \psi,\qquad \text{ in }(-(2\upsilon)^2, 0)\times \mathbb{R}^{3d},\\
g(-(2\upsilon)^2, \cdot )&= 0,  \qquad\quad\text{ in }\mathbb{R}^{3d}.
\end{cases}
\end{equation}

From Lemma \ref{F} we know 
\begin{equation}\label{g1}
  \big(| \dz g |^2\big) _{Q_{\upsilon}} ^{1/2}\leq N \sum\limits_{k=0}^{\infty} 2^{-2k} \big(| f |^2 \big) _{Q_{\upsilon,2^{(k+1)\upsilon}}}^{1/2}.
\end{equation}

Besides, by H\"{o}der inequality we have
\begin{equation}\label{g2}
\begin{aligned}
 \big(| \dz g |^2\big)_{Q_{1}}^{1/2}&\leq N \upsilon^{\frac{2+9d}{2}} \big(| \dz g |^2\big)^{1/2} _{Q_{\upsilon}}\\
 & \leq N \upsilon^{\frac{2+9d}{2}}\sum\limits_{k=0}^{\infty} 2^{-2k} \big(| f |^2 \big) _{Q_{\upsilon,2^{(k+1)}\upsilon}}^{1/2}.
 \end{aligned}
\end{equation}

Next, we consider the equation satisfied by $h = u - g$
\[ P_0 h = f (1 - \psi).\]
Note that $\psi = 1$ on $( -\upsilon^2, 0) \times B_{\upsilon} \times B_{\upsilon^3}$. Combining Proposition~\ref{sharp} with estimate \eqref{g1}, we derive
\begin{equation}
\begin{aligned}
&\big(|\dz h - (\dz h)_{Q_1}|^2\big)^{1/2}_{Q_1} \\
 \leq & N\upsilon^{-1} \big( |\dz h|^2\big)^{1/2}_{Q_{\upsilon}} \\
 \leq & N\upsilon^{-1} \big( |\dz u|^2\big)^{1/2}_{Q_{\upsilon}}  +N\upsilon^{-1} \big( |\dz g|^2\big)^{1/2}_{Q_{\upsilon}}\\
 \leq & N \upsilon^{-1}\big( |\dz u|^2\big)^{1/2}_{Q_{\upsilon}}+ N\upsilon^{-1}\sum\limits_{k=0}^{\infty} 2^{-2k}\big(| f |^2 \big) _{Q_{\upsilon,2^{(k+1)\upsilon}}}^{1/2}.
 \end{aligned}
 \end{equation}
 
The combination of this inequality with \eqref{g2} directly yields \eqref{sharp:dzu}.
 
The term $I_2$ is handled similarly through another application of Lemma~\ref{F}, yielding
\begin{equation}\label{g3}
\big( | \nabla^2_x g |^2\big)_{Q_{\upsilon}} ^{1/2}\leq N \sum\limits_{k=0}^{\infty} 2^{-k^2/8}\big(|f|^2\big) ^{1/2}_{  Q_{\upsilon,2^{k+1}\upsilon}}.
\end{equation}

Then, we have
\begin{equation}\label{g4}
\begin{aligned}
\big( | \nabla^2_x g |^2\big) _{Q_{1}} ^{1/2}&\leq N \upsilon^{\frac{2+9d}{2}}\big( | \nabla^2_x g |^2\big)  _{Q_{\upsilon}} ^{1/2}\\
&\leq N \upsilon^{\frac{2+9d}{2}} \sum\limits_{k=0}^{\infty} 2^{-k^2/8} \big(|f|^2\big) ^{1/2}_{  Q_{\upsilon,2^{k+1}\upsilon}}.
\end{aligned}
\end{equation}

By Lemma \ref{sharp}, we have the estimate for $\nabla^2_x h$
\begin{equation}
\begin{aligned}
&\big(|\nabla^2_x h - (\nabla^2_x h)_{Q_1}|^2\big)^{1/2}_{Q_1} \\ \leq & N \upsilon^{-1} \big(|\nabla^2_x h |^2\big)^{1/2}_{Q_{\upsilon}} +N \upsilon^{-1} \sum\limits_{k=0}^{\infty} 2^{-3k}\big( |\dz h|^2\big)^{1/2}_{Q_{\upsilon , 2^{k}\upsilon }}\\
\leq &N \upsilon^{-1} \big( |\nabla^2_x u |^2\big)^{1/2}_{Q_{\upsilon}} + N \upsilon^{-1} \sum\limits_{k=0}^{\infty} 2^{-3k} \big( |\dz u|^2\big)^{1/2}_{Q_{\upsilon , 2^{k}\upsilon }}\\
 &+ N \upsilon^{-1} \big( |\nabla^2_x g |^2\big)^{1/2}_{Q_{\upsilon}} + N \upsilon^{-1} \sum\limits_{k=0}^{\infty} 2^{-3k} \big( |\dz g|^2\big)^{1/2}_{Q_{\upsilon , 2^{k}\upsilon }}\\
 \leq & N \upsilon^{-1}  \big( |\nabla^2_x u |^2\big)^{1/2}_{Q_{\upsilon}} + N \upsilon^{-1} \sum\limits_{k=0}^{\infty} 2^{-3k}  \big( |\dz u|^2\big)^{1/2}_{Q_{\upsilon , 2^{k}\upsilon }}\\
&+N  \sum\limits_{k=0}^{\infty} 2^{-k^2/8} \big(|f|^2\big) ^{1/2}_{  Q_{\upsilon,2^{k}\upsilon}} + N \upsilon^{-1}\sum\limits_{k=0}^{\infty} 2^{-k} \big(|f|^2\big)^{1/2}_{Q_{\upsilon,2^{k}\upsilon}}.
 \end{aligned}
 \end{equation}
Combining this inequality with \eqref{g4} yields \eqref{sharp:dzu}.
\end{pf}

\vskip .1in
With pointwise estimates for the sharp functions $\dz u$ and $\nabla_x^2 u$ now established, we shall apply the Hardy-Littlewood maximal theorem and Fefferman-Stein sharp function theory to derive global $L^p$ estimates for these derivatives.

\begin{proposition}\label{zx}
For any $p \in(2, \infty) ,$ \T . Suppose $u \in S^{p}(\mathbb{R}^{1+3d}_T)$, then we have
\[ \Vert \nabla^2_x u \Vert_{L^{p}( \mathbb{R}_T^{1+3d})} +  \Vert \dz u \Vert_{L^{p}( \mathbb{R}_T^{1+3d})} \leq N(d,\delta,p) \Vert P_0 u \Vert_{L^{p}( \mathbb{R}_T^{1+3d})}.\]
\end{proposition}
\begin{pf}
By Lemma \ref{sharp} we conclude that
\begin{equation}
  \begin{aligned}
&\big(\dz u \big)^{\sharp}_{T}(X)\\
\leq & N  \upsilon^{-1}\mathcal {M}_T^{1/2}|\dz u |^2 (X)+ N  \upsilon^{\frac{2+9d}{2}}\sum\limits_{k=0}^{\infty} 2^{-2k}\mathcal {M}_{2^k,T}^{1/2}|f |^2(X),\\
&\big(\nabla^2_x u\big)^{\sharp}_{T}(X)\\
\leq & N  \upsilon^{-1}\mathcal {M}_T^{1/2}|\nabla^2_x u |^2(X) + N  \upsilon^{-1}\sum\limits_{k=0}^{\infty} 2^{-3k}\mathcal {M}_{2^k,T}^{1/2}|\dz u|^2(X) \\
&+N \upsilon^{\frac{2+9d}{2}}\sum\limits_{k=0}^{\infty} 2^{-k}\mathcal {M}_{2^k,T}^{1/2}|f |^2(X).
\end{aligned}
\end{equation}

Applying the Hardy-Littlewood theorem and the Fefferman-Stein theorem, we obtain
\begin{equation}\label{dzu}
\begin{aligned}
&\Vert \dz u \Vert_{L^{p}( \mathbb{R}_T^{1+3d})} \\
\leq & N \upsilon^{-1} \Vert \dz u \Vert_{L^{p}( \mathbb{R}_T^{1+3d})}+ N \upsilon^{\frac{2+9d}{2}} \Vert f \Vert_{L^{p}( \mathbb{R}_T^{1+3d})},\\
&\Vert \nabla^2_x u \Vert_{L^{p}( \mathbb{R}_T^{1+3d})}\\
\leq & N \upsilon^{-1} \Vert \nabla^2_x u \Vert_{L^{p}( \mathbb{R}_T^{1+3d})}+N \upsilon^{-1} \Vert \dz u \Vert_{L^{p}( \mathbb{R}_T^{1+3d})}+ N \upsilon^{\frac{2+9d}{2}} \Vert f \Vert_{L^{p}( \mathbb{R}_T^{1+3d})}.
\end{aligned}
\end{equation}

Let $\upsilon = 2N+2$ in \eqref{dzu}, we get the estimate of $\dz u$ and $\nabla_x^2 u$. The Proposition has been proved.
\end{pf}

\begin{lemma}\label{samgon}
Under the assumptions of Proposition \ref{zx}, for any $\lambda \geq 0$, we have
\begin{equation}\label{u}
\lambda \Vert u \Vert_{L^{p}( \mathbb{R}_T^{1+3d})} \leq N( d,\delta,p)\Vert P_0 u + \lambda u \Vert_{L^{p}( \mathbb{R}_T^{1+3d})}.
\end{equation}
\end{lemma} 
\begin{pf}
  Denote
  \[ \hat{x} = (x_1, \cdots,x_{d+1}),  \hat{y} = (y_1, \cdots,y_{d+1}), \hat{z} = (z_1, \cdots,z_{d+1}).\]
$$ \hat{P}_0(\hat{X}) = \partial_t - \sum_{i=1}^{d+1} x_i \partial_{y_i} - \sum_{i=1}^{d+1} y_i \partial_{z_i} - \sum_{i,j=1}^{d} a^{ij}(t)\partial_{x_{i}x_{j}} - \partial_{x_{d+1}x_{d+1}}.$$

Let $\psi \in C^{\infty}_{0}(\mathbb{R})$ and $\psi\neq0$. Set
$$ \hat{u} (\hat{X}) = u(X) \psi(x_{d+1})cos(\lambda^{1/2}x_{d+1}).$$

Then by direct calculation, we have 
\begin{equation}\label{miss}
\begin{aligned}
\partial_{x_{d+1}x_{d+1}} \hat{u} (\hat{X}) &= u(X) \psi^{\prime \prime}(x_{d+1})cos(\lambda^{1/2}x_{d+1})- \lambda u(X) \psi(x_{d+1})cos(\lambda^{1/2}x_{d+1})
\\
&-2\lambda^{1/2}u(X) \psi^{\prime}(x_{d+1})sin(\lambda^{1/2}x_{d+1}).
\end{aligned}
\end{equation}


Furthermore, we conclude that
\begin{equation}\label{hatP}
\begin{aligned}
\hat{P}_0 \hat{u} (\hat{X}) = &P_0 u(X) \psi(x_{d+1})cos(\lambda^{1/2}x_{d+1})- u(X)\psi^{\prime \prime}(x_{d+1})cos(\lambda^{1/2}x_{d+1})  \\
& + \lambda \hat{u}(\hat{X})-2\lambda^{1/2}u(X) \psi^{\prime}(x_{d+1})sin(\lambda^{1/2}x_{d+1}).
\end{aligned}
\end{equation}

Note for all $p>0$ and $\lambda >1$, we have
\[ \int_{\mathbb{R}}|\psi(t)cos(\lambda^{1/2}x_{d+1})|^p \mathrm{d}t \geq N(p)>0.\]

Then combined with \eqref{hatP}, we have
$$\lambda \Vert u \Vert_{L^{p}(\mathbb{R}_T^{1+3d})} \leq N  \Vert \partial_{x_{d+1}x_{d+1}} \hat{u} (\hat{X}) \Vert_{L^{p}(\mathbb{R}_T^{1+3d})} + N(1+\lambda^{1/2})\Vert u \Vert_{L^{p}(\mathbb{R}_T^{1+3d})}.$$

With the help of  \eqref{hatP} and by Proposition \ref{zx},
\begin{equation}
\begin{aligned}
 &\Vert\partial_{x_{d+1}x_{d+1}}\hat{u} \Vert_{L^{p}(\mathbb{R}_T^{1+3d})}
 \leq  N \Vert \hat{P}_0 \hat{u} \Vert_{L^{p}(\mathbb{R}_T^{1+3d})}\\
  \leq & N \Vert P_0 u +\lambda u \Vert_{L^{p}(\mathbb{R}_T^{1+3d})} + N(1+\lambda^{1/2})\Vert u \Vert_{L^{p}(\mathbb{R}_T^{1+3d})}.
  \end{aligned}
  \end{equation}
 
 That gives
 \begin{equation}\label{lambda_u}
 \lambda \Vert u \Vert_{L^{p}(\mathbb{R}_T^{1+3d})} \leq N \Vert P_0 u +\lambda u \Vert_{L^{p}(\mathbb{R}_T^{1+3d})} + N(1+\lambda^{1/2})\Vert u \Vert_{L^{p}(\mathbb{R}_T^{1+3d})}.
 \end{equation}

If we choose $\lambda$ lager enough such that $\lambda \geq \lambda_0$, where $\lambda_0 = 16 N^2 +1$. Then $\lambda - N(1+\lambda^{1/2}) > \lambda/2$. By \eqref{lambda_u}, we get
 $$ \lambda \Vert u \Vert_{L^{p}(\mathbb{R}_T^{1+3d})} \leq N \Vert P_0 u +\lambda u \Vert_{L^{p}(\mathbb{R}_T^{1+3d})}. $$
 
 Using scaling we can also get the desired estimate for $0 <\lambda < \lambda_0$.
 
\end{pf}

\vskip.1in

Following the approach of Lemma~\ref{local}, we can also establish localized $L^p$ estimates for the solution $u$.
\begin{lemma}\label{localp}
  Let $\lambda \geq 0$, $0< r_1 < r_2,0 <  R_1 <  R_2$, $p\leq 2$, assume $u \in S^{p}_{loc} (\mathbb{R}^{1+3d}_{0})$ . Denote $f = P_0 u + \lambda u$, then there exist a constant $N=N(d,\delta)$ such that the following local estimates hold.
  \begin{equation}
    \begin{aligned}
    (\romannumeral 1) & ( r_2 - r_1  )^{-1} \Vert \nabla_x u \Vert_{L^{p}(Q_{r_1, R_1})} +  \Vert \nabla_x ^2 u \Vert_{L^{p}( Q_{r_1, R_1})}\\
    \leq & N(d,\delta,p)\Big(\Vert f \Vert_{L^{p}(Q_{r_2, R_2})}\\
    & + (( r_2 - r_1 )^{-2} + r_2 (R_2 - R_1)^{-3} +  R_2 (R_2 - R_1)^{-5})\Vert u \Vert_{L^{p}(Q_{r_2, R_2})}\Big).
  \end{aligned}
\end{equation} 
  \qquad$(\romannumeral 2)$ Denote $C_r = ( -r^2 ,0)\times B_r \times B_{r^3}\times \mathbb{R}^{d}.$ Then we have
  \begin{equation}
    \begin{aligned}
      &( r_2 - r_1  ) ^{-1} \Vert \nabla_x u \Vert_{L^{p}(C_{r})} +  \Vert \nabla_x ^2 u \Vert_{L^{p}(C_{r})}\\
       \leq & N(d,\delta,p)\Big(\Vert f \Vert_{L^{p}(C_{r})} +  (r_2 - r_1)^{-2} \Vert u \Vert_{L^{p}(C_{r})}\Big).
    \end{aligned}
  \end{equation}
\end{lemma}

We generalize Lemma~\ref{dense2} to establish the corresponding $L^p$ density property.
\begin{lemma}\label{densep}
For any $\lambda \ge 0$ and $p >1$,
the set $(P_0 + \lambda) C^{\infty}_0 (\mathbb{R}^{1+3d})$
is dense in $L^p (\mathbb{R}^{1+3d})$.
\end{lemma}

\vskip .1in
\no{\bf Proof of Theorem \ref{Thm3.1}.}\quad
First we consider the case $p>2$. Combined Proposition \ref{zx} with Lemma \ref{samgon}, we conclude that the estimates for $\lambda u$, $\nabla_x^2$ and $\dz u$ hold. Throughout the proof, we assume that $N = N (d, p)$. By interpolation inequality we have 
  $$ \lambda^{1/2}\|\nabla_x u\|_{L^{p}({\mathbb{R}^{1+3d}_{T}} )}\leq \lambda \| u\|_{L^{p}({\mathbb{R}^{1+3d}_{T}} )} + N  \| \nabla_x^2u\|_{L^{p}({\mathbb{R}^{1+3d}_{T}} )}.$$
  Note that $\{a^{ij}\}$ is a time-dependent matrix, by Theorem 1.1 of \cite{chen2019propagation}, one has
  $$ \|\dy u\|_{L^{p}({\mathbb{R}^{1+3d}_{T}} )}\leq N \| P_0 u + \lambda u\|_{L^{p}({\mathbb{R}^{1+3d}_{T}} )}.$$
  
  As for $\nabla_x(-\Delta_y)^{1/6}$, thanks to Appendix \ref{Hormander}, we have 
  $$ \|\nabla_x(-\Delta_y)^{1/6}\|_{L^{p}({\mathbb{R}^{1+3d}_{T}} )}\leq N  \| \nabla_x^2u\|_{L^{p}({\mathbb{R}^{1+3d}_{T}} )}+ \| \dy u\|_{L^{p}({\mathbb{R}^{1+3d}_{T}} )}.$$
  Furthermore combined Lemma \ref{densep} with the prior estimate, we obtain the existence and uniqueness of the equation \eqref{maineq}.
  
  For the range $p \in (1,2)$, the result follows by a duality argument analogous to Theorem 5.1 in \cite{dong2022global}. Here we just omit it.
\endproof

Building on Theorem~\ref{Thm3.1}, we now extend the preceding lemmas and propositions to the full range $p \in (1,\infty)$.

\begin{lemma}\label{localpp}
The assertions of Lemma \ref{localp} and Lemma \ref{F}- Lemma \ref{x2} hold for any $p \in (1,\infty)$.
\end{lemma}
\begin{pf}
By adapting Lemma~\ref{localp} through substitution of Theorem~\ref{Thm:L2} with Theorem~\ref{Thm3.1}, we similarly obtain localized $L^p$ estimates.
\end{pf}

\begin{lemma}
Lemma \ref{F}- Lemma \ref{x2} hold for $p>1$.
%
%
%
\end{lemma} 
\begin{pf}
The proofs follow the same arguments as in the corresponding lemmas, requiring only the substitution of Theorem~\ref{Thm:L2} by Theorem~\ref{Thm3.1} and Lemma~\ref{local} by Lemma~\ref{localpp}.
\end{pf}
Furthermore, we can also establish a generalization of Proposition~\ref{sharp} valid for all $p \in (1,\infty)$.
\begin{proposition}\label{psharp}
Let $p>1$, $r>0$, $\upsilon\geq 2$, $T\in ( -\infty, \infty]$, $X_0 \in \overline{\mathbb{R}^{1+3d}_T}$. Suppose $u \in S^p (\mathbb{R}^{1+3d}_T)$. Assume $P_0 u =f$ in $\mathbb{R}^{1+3d}_T$. Then there exits a constant $N =N(d,\delta,p)$, so that
\begin{align*}
    (\romannumeral 1)\,\,\,&\big(|\dz u - (\dz u) _{Q_r(X_0)}|^p \big) ^{1/p}_{Q_r(X_0)}\\
    \leq & N \upsilon^{-1} \big( |\dz u|^p\big)^{1/p}_{Q_{\upsilon r}(X_0)} + N \upsilon^{\frac{2+9d}{2}} \sum\limits_{k=0}^{\infty} 2^{-2k} \big(|f|^p\big)_{Q_{\upsilon r,2^k\upsilon r}(X_0)}^{1/p},\\
    (\romannumeral2)\,\,\,&\big(|\nabla_x^2 u - \big(\nabla_x^2 u\big)_{Q_r(X_0)}|^p\big)^{1/p}_{Q_r(X_0)} \\
   \leq & N \upsilon^{-1} ( |\nabla_x^2 u|^p)^{1/p}_{Q_{\upsilon r}(X_0)} + N \upsilon^{-1} \sum\limits_{k=0}^{\infty} 2^{-3k} \big( |\dz u|^p\big)^{1/p}_{Q_{\upsilon r, 2^k\upsilon r}(X_0)}\\
   & + N\upsilon^{\frac{2+9d}{2}}  \sum\limits_{k=0}^{\infty} 2^{-k} \big(|f|^p\big)  _{Q_{\upsilon r,2^k\upsilon r}(X_0)}^{1/p}.
\end{align*}
\end{proposition} 

\vskip .2in
\section{The proof of the main result}\hspace*{\parindent}

This section addresses coefficients $a^{ij}(X)$ satisfying assumption $\mathbf{[A_2]}$. Employing the frozen coefficient method and building on Section 3's results, we establish sharp function estimates for $\nabla_x^2u$ under VMO conditions on $a^{ij}$.
\begin{lemma}\label{adx}
Let $\theta_0 > 0$, $\upsilon\geq 2$, $ \alpha \in (1, 5/3)$, $q\in (2,\infty)$,\T. Assume $R_0$ be the constant of $\mathbf{[A_2]}$. Suppose $u \in S^{q}(\mathbb{R}^{1+3d}_{T})$, then there exist a constant $N = N (d,\delta, p)$ and a sequence $\{a_k, k\geq0\}$ and 
\[ \sum\limits_{k=0}^{\infty} a_k \leq N.\]
For any $X_0\in \overline{\mathbb{R}^{1+3d}_T} $ , $r \in (0, R_0/(4\upsilon))$, we have
\begin{equation}\label{d2x}
\begin{aligned}
&\left(|\nabla_x^2 u - (\nabla_x^2 u ) _{Q_r(X_0)}|^q\right)^{1/q}_{Q_r(X_0)} \\
 \leq & N  \upsilon^{-1} (|\nabla_x^2 u|^q)^{1/q}_{Q_{\upsilon r}(X_0)} +  N\upsilon^{-1} \sum\limits_{k=0}^{\infty} 2^{-3k} \left( |\dz u|^q\right)^{1/q}_{Q_{\upsilon r, 2^k\upsilon r}(X_0)}\\
& + N\upsilon^{\frac{2+9d}{2}}  \sum\limits_{k=0}^{\infty} 2^{-k} ( f^q ) _{Q_{2\upsilon r,2^{k}\upsilon r}(X_0)}^{1/q}+ N \upsilon^{\frac{2+9d}{2}}\theta_0 ^{( \alpha-1)/(q\alpha)}  \sum\limits_{k=0}^{\infty} a_k\left(|\nabla^2_xu|^{q\alpha}\right) _{Q_{2\upsilon r,2^{k }\upsilon r} }^{1/(q\alpha)}.
\end{aligned}
\end{equation}
\end{lemma} 

Following Lemma 7.2 of \cite{dong2022global}, we can reformulate assumption $\mathbf{[A_2]}$ as:
\begin{lemma}\label{a}
Let $\theta_0 > 0$, $R_0$ be the constants in $\mathbf{[A_2]}$, for $r \in ( 0, R_0/2)$, $c\geq 1$, we have
\begin{equation}\label{Id}
  I:=\fint _{Q_{r,cr}}|a( t,x,y,z)-(a(t, \cdot,\cdot,\cdot))_{B_r\times B_{r^3}\times B_{r^5}}|dX \leq N c^5 \theta_0.
\end{equation}
\end{lemma}

 \vskip .1in
 \no{\bf Proof of Theorem \ref{T1}.\quad}
At first we consider the situation that $|\vec{b}|=c=0$. Suppose that $u \in S^{p}(\mathbb{R}^{1+3d}_{T})$.

Let $1 < q < p$ and $t_0 \in \mathbb{R}$, consider solutions $u$ vanishing outside the cylindrical domain $(t_0 - (R_0R_1)^2, t_0) \times \mathbb{R}^{3d}$.

When $4\upsilon r \geq R_0$, an application of Hölder's inequality yields, for all $X \in \overline{\mathbb{R}_{T}^{1+3d}}$
 \begin{equation}
  \begin{aligned}
    \big(|\nabla_x^2 u - ( \nabla_x^2 u) _{Q_r(X)}|^q\big)^{1/q}_{Q_r(X)} & \leq 2  (|\nabla_x^2 u |^q)^{1/q}_{Q_r(X)} \\
  &\leq 2(\chi_{(t_0 - (R_0R_1)^2,t_0)})_{Q_r(X)}^{1/q\alpha_1}(|\nabla_x^2 u |^{2\alpha})^{1/q\alpha}_{Q_r(X)}\\
  & \leq 2 (R_1R_0r^{-1})^{2/q\alpha_1} \mathcal {M}_{T}^{1/\zk q \alpha\yk}|\nabla_x^2u|^{q\alpha}(X)\\
  & \leq N \upsilon^{2/q\alpha_1}R_1^{2/q\alpha_1}  \mathcal {M}_{T}^{1/(q \alpha)}|\nabla_x^2u|^{q\alpha}(X).
  \end{aligned}
 \end{equation}
 
 For the complementary case $4\upsilon r < R_0$, Lemma~\ref{adx} applies. Combining these two cases yields:
 \begin{equation}
 \begin{aligned}
 \zk \nabla^2_x u \yk ^{\sharp}_{T}(X)& \leq N \upsilon^{-1} \mathcal {M}_{T}^{1/q}|\nabla_x^2u|^{q}(X)
 +  N \upsilon^{2/q\alpha_1}R_1^{2/q\alpha_1}  \mathcal {M}_{T}^{1/(q \alpha)}|\nabla_x^2u|^{q\alpha}(X)\\
 &+ N \upsilon^{\frac{2+9d}{2}}\theta_0 ^{(\alpha-1)/(q\alpha)}  \sum\limits_{k=0}^{\infty} a_k  \mathcal {M}_{2^k,T}^{1/(q\alpha)}|\nabla_x^2u|^{q\alpha}(X)\\
 &+ N\upsilon^{-1} \sum\limits_{k=0}^{\infty} 2^{-3k} \mathcal {M}_{2^k,T}^{1/q}|\dz u|^{2}(X)\\
 &+ N\upsilon^{\frac{2+9d}{2}}  \sum\limits_{k=0}^{\infty} 2^{-k}\mathcal {M}_{2^k,T}^{1/q}|Pu|^{q}(X).
 \end{aligned}
 \end{equation}
 
Taking the $L^p$ norm on both sides and applying Minkowski's inequality, we derive:
\begin{equation}\label{fd2x}
\begin{aligned}
  &\Vert \nabla_x^2 u \Vert _{L^{p}(\mathbb{R}^{1+3d}_T)}\\
  \leq & N \upsilon^{-1}\Vert \nabla_x^2 u \Vert _{L^{p}(\mathbb{R}^{1+3d}_T)}+N \upsilon^{1/\alpha_1}R_1^{1/\alpha_1}\Vert \nabla_x^2 u \Vert _{L^{p}(\mathbb{R}^{1+3d}_T)}\\
  &+N\upsilon^{\frac{2+9d}{2}}\theta_0 ^{\zk \alpha-1\yk/\zk2\alpha\yk}\Vert \nabla_x^2 u \Vert _{L^{p}(\mathbb{R}^{1+3d}_T)}+ N \upsilon^{-1} \Vert \dz u \Vert _{L^{p}(\mathbb{R}^{1+3d}_T)}\\
  & + N \upsilon^{\frac{2+9d}{2}}  \Vert P u \Vert _{L^{p}(\mathbb{R}^{1+3d}_T)}.
  \end{aligned}
\end{equation}

We now estimate $\dz u$. Observing that $u$ satisfies:
\[ \partial_t u - x\cdot\nabla_y u -y\cdot\nabla_z - \nabla^2_x u = Pu + (a^{ij}-\delta^{ij}) \partial_{x_i}u\partial_{x_j}u,\]

By Theorem \ref{Thm3.1}, we obtain
\begin{equation}\label{fdz}
\Vert \dz u \Vert _{L^{p}(\mathbb{R}^{1+3d}_T)} \leq N  \Vert P u \Vert _{L^{p}(\mathbb{R}^{1+3d}_T)} + N\Vert \nabla_x^2 u \Vert _{L^{p}(\mathbb{R}^{1+3d}_T)}.
\end{equation}

Back to \eqref{fd2x}, we find 
\begin{equation}\label{ffd2x}
\begin{aligned}
  &\Vert \nabla_x^2 u \Vert _{L^{p}(\mathbb{R}^{1+3d}_T)}\\
  \leq & N \upsilon^{-1}\Vert \nabla_x^2 u \Vert _{L^{p}(\mathbb{R}^{1+3d}_T)}+N \upsilon^{1/\alpha_1}R_1^{1/\alpha_1}\Vert \nabla_x^2 u \Vert _{L^{p}(\mathbb{R}^{1+3d}_T)}\\
  &+N\upsilon^{\frac{11}{2}}\theta_0 ^{\zk \alpha-1\yk/\zk2\alpha\yk}\Vert \nabla_x^2 u \Vert _{L^{p}(\mathbb{R}^{1+3d}_T)}
   + N \zk \upsilon^{\frac{11}{2}} + \upsilon^{-1}\yk  \Vert P u \Vert _{L^{p}(\mathbb{R}^{1+3d}_T)}.
  \end{aligned}
\end{equation}

Choose $\upsilon = 2 + 4N$, $\theta_0 >0 $, $R_1 > 0$ small enough such that
\[ N \upsilon^{1/\alpha_1}R_1^{1/\alpha_1} \leq 1/4, \qquad N\upsilon^{\frac{11}{2}}\theta_0 ^{\zk \alpha-1\yk/\zk2\alpha\yk}\leq 1/4.\]

By eliminating the term $\nabla_x^2 u$ from the right-hand side, we obtain the estimate of $ \nabla_x^2 u$. Then, we can derive the desired estimate for $\| \dz u \|_{L^{p}(\mathbb{R}^{1+3d}_T)}$ from \eqref{fdz}.

According to Theorem 1 of \cite{chen2019propagation}, we have
\begin{equation}\label{fdy}
\begin{aligned}
&\Vert \dy u \Vert _{L^{p}(\mathbb{R}^{1+3d}_T)}
 \leq N  \Vert P u \Vert _{L^{p}(\mathbb{R}^{1+3d}_T)} + N\Vert \nabla_x^2 u \Vert _{L^{p}(\mathbb{R}^{1+3d}_T)}\\
 \leq& N \Vert P u \Vert _{L^{p}(\mathbb{R}^{1+3d}_T)}.
\end{aligned}
\end{equation}

The estimate for $\nabla_x(-\Delta_y)^{1/6}$ follows from the interpolation inequality in Appendix~\ref{Hormander}:
\begin{equation}\label{fxdy}
\begin{aligned}
&\Vert \nabla_x(-\Delta_y)^{1/6} u \Vert _{L^{p}(\mathbb{R}^{1+3d}_T)}\\
 \leq &N  \Vert \dy u \Vert _{L^{p}(\mathbb{R}^{1+3d}_T)} + N\Vert \nabla_x^2 u \Vert _{L^{p}(\mathbb{R}^{1+3d}_T)}
 \leq N \Vert P u \Vert _{L^{p}(\mathbb{R}^{1+3d}_T)}.
 \end{aligned}
\end{equation}

For general solutions $u$, the result follows by applying a cutoff function and adapting the method of Theorem 2.6 in \cite{dong2022global}.

\vskip .3in

\appendix
\section{}
\renewcommand{\appendixname}{\appendixname~\Alph{section}}
\setcounter{equation}{0}

\begin{lemma}\label{my1}
  Assume $u \in C^{\infty}(Q_1 )$ and $P_0 u =0 $ in $Q_1$. Suppose that $\nabla^2_x u = \nabla_z u = 0$. Besides for $i=1,\cdots,d$, we have 
  \begin{equation}\label{xyu}
        \int_{Q_1} u = \int_{Q_1} x_i u = \int_{Q_1} y_iu = 0.
        \end{equation}
        Besides for $1\leq i<j\leq d$,
        \begin{equation}\label{xyu=0}
        \int_{Q_1} x_iy_j u  = 0.
        \end{equation}
        Then we get that $u\equiv0$ in $Q_1$ .
\end{lemma}
\begin{pf}
        Note that $P_0 u =0 $  in $Q_1$ and $\partial^2_x u = \partial_z u = 0$, imply
          \[ (\partial_t - x \cdot \nabla_y ) u =0.\]
          
          Set $v( s) = u ( s, x, (t-s) x + y)$, we have 
          \[ \frac{\mathrm{d}v(s) }{\mathrm{d}s} \equiv 0.\]
          So we get $u(t, x, y) = g(x, tx + y) = :u(0, x, tx+y)$.

            Next we use  $\nabla^2_x u =0 $ to get the representation of $u$.
          Since
          \begin{align*}
            0 &= \partial_{x_ix_j} u (t, x, y)\\
            &= \partial_{i,j}g(x, tx + y)+ t \partial_{i,j+d}g(x, tx + y) +t \partial_{i+d, j}g(x, tx + y) + t^2 \partial_{i+d,j+d}g(x, tx + y),
          \end{align*}
         where $\partial_{i} g$ is the derivative of the $i-th$ component of $g$.
         
         Let $t\rightarrow0$, we have
         \[ \partial_{ij}g(x, y) = 0. \]
         Then we can find $b_0, b_i$ such that $g(x, y)= b_0( y) + \sum_{i=1}^{d}b_i(y)x_i. $
         That shows
         \begin{equation}\label{uh}
         u (t, x, y)= b_0( y+tx) +\sum_{i=1}^{d}b_i(y+tx)x_i.
         \end{equation}
         
        By $\nabla^2_x u = 0 $ again, 
         \[ t^2 \partial_{kl}b_0(tx+y) + t \partial_k b_{l}(tx+y)+ t \partial_l b_{k}(tx+y) +t^2 \sum_{i=1}^{d}\partial_{kl}b_{i}(tx+y)x_i \equiv 0. \]
         
        Let $x = 0$, we get
        \begin{equation}
        \begin{cases}
        \partial_{kl}b_0= 0,\\
        \partial_k b_{l} + \partial_l b_{k} = 0.
        \end{cases}
        \end{equation}
        Then we get 
        \begin{equation}
        \begin{cases}
        b_0(y) = c_0 + \sum_{i=1}^{d}c_i y_i,\\
        b_l (y) = h_l + \sum_{i=1}^{d}h_{li}y_i.
        \end{cases}
        \end{equation}
        where $h_{lj}+h_{jl} =0$.
        Back to \eqref{uh}, we conclude that
         \begin{equation}
         u(t, x, y) = c_0 + \sum_{i=1}^{d}c_i (y_i +tx_i) + \sum_{i=1}^{d}h_i x_i +  \sum_{1\leq i<j\leq d}h_{ij}(x_iy_j - x_jy_i).
         \end{equation}
        According to \eqref{xyu=0}, for $1\leq i<j\leq d$, we obtain
        $$ \int_{Q_1} x_i y_j u = h_{ij}\int_{Q_1} x_i^2y_j^2 =0,$$
        that implies $h_{ij} =0$.
        By \eqref{xyu}, for $1\leq i \leq d$, we have
        $$ \int_{Q_1} y_i u = c_i \int_{Q_1}y_i^2 = 0,$$
        so $c_i = 0.$
        Similarly, 
        $$ \int_{Q_1} u = \int_{Q_1} c_0 = \int_{Q_1} x_i u = h_i \int_{Q_1} x_i^2 =0,$$
        we obtain that $c_0 = h_i= 0$. Therefore $u\equiv0$.   
\end{pf}

\begin{lemma}\label{Hormander}
For any $p \in (1,\infty)$ and functions $u(x,y)$ defined on $\mathbb{R}^{2d}$, the following interpolation inequality holds
$$ \|\nabla_x(-\Delta_y)^{1/6}u\|_{L^p(\mathbb{R}^{2d})}\leq N(d,p) \Big(\|\nabla^2_xu\|_{L^p(\mathbb{R}^{2d})}+\|(-\Delta_y)^{1/3}u\|_{L^p(\mathbb{R}^{2d})} \Big).$$
\end{lemma}
\begin{pf}
Denote $\cF{h}(\xi,\eta)$ as the Fourier transform of $h(x,y)$. Then 
$$ \cF{\nabla_x(-\Delta_y)^{1/6}u}= \xi |\eta|^{1/3} \cF{u} = \frac{\xi |\eta|^{1/3}}{|\xi|^2+|\eta|^{2/3}}(\cF\nabla^2_xu+\cF(-\Delta_y)^{1/3}u).$$
Set $m(\xi,\eta)= \frac{\xi |\eta|^{1/3}}{|\xi|^2+|\eta|^{2/3}}$, then for any $k>0$, one has
$$ m(k\xi,k^2\eta)=m(\xi,\eta).$$
Note that $m$ is a bounded function on $\mathbb{R}^{2d}$, therefore by Corollary 6.2.5 of \cite{grafakos2008classical},  $m$ is a Marcinkiewicz Multiplier on $\mathbb{R}^{2d}$. Thus we conclude that 
$$ \|\nabla_x(-\Delta_y)^{1/6}u\|_{L^p(\mathbb{R}^{2d})}\leq N(d,p) \Big(\|\nabla^2_xu\|_{L^p(\mathbb{R}^{2d})}+\|(-\Delta_y)^{1/3}u\|_{L^p(\mathbb{R}^{2d})} \Big).$$
\end{pf}


\section*{Acknowledgments}
The author is deeply grateful to Professor Zhang for the fundamental insight of employing proof by contradiction. Special gratitude is extended to Professor Dong for his invaluable suggestions on the Poincaré inequality.

\vskip .3in
\bibliographystyle{plain}
\bibliography{ref}

\end{document}